\documentclass{amsart}

\usepackage{amsmath, amscd, amssymb}
\usepackage{graphicx, psfrag}
\usepackage{color}
\newcommand{\cal}{\mathcal}

\newcommand{\bC}{{\Bbb C}}
\newcommand{\bE}{{\Bbb E}}
\newcommand{\bF}{{\Bbb F}}
\newcommand{\bL}{{\Bbb L}}

\newcommand{\bP}{{\Bbb P}}
\newcommand{\bQ}{{\Bbb Q}}
\newcommand{\bR}{{\Bbb R}}
\newcommand{\bZ}{{\Bbb Z}}

\newcommand{\cD}{{\cal D}}
\newcommand{\cE}{{\cal E}}
\newcommand{\cF}{{\cal F}}

\newcommand{\cM}{{\cal M}}
\newcommand{\cN}{{\cal N}}
\newcommand{\cO}{{\cal O}}
\newcommand{\cP}{{\cal P}}

\newcommand{\cU}{{\cal U}}

\newcommand{\cW}{{\cal W}}

\newcommand{\vir}{ {\mathrm{vir}} }

\newcommand{\tpi  }{\tilde{\pi}  }

\newcommand{\tF}{\tilde{F}}

\newcommand{\tV}{\tilde{V}}

\newcommand{\tf}{\tilde{f}}
\newcommand{\tz}{\tilde{z}}

\newcommand{\MgY}{\Mbar_{g,0}(Y; D^1,\ldots,D^k \mid \beta; \mu^1,\ldots, \mu^k)}
\newcommand{\GgY}{ G_{g,0}(Y;D^1,\ldots, D^k\mid \beta;\mu^1,\ldots,\mu^k) }
\newcommand{\MY}{\Mbar^\bullet_\chi(Y; D^1,\ldots,D^k \mid \beta;\mu^1,
\cdots,\mu^k)}
\newcommand{\GY}{G^\bullet_\chi(Y; D^1,\ldots,D^k \mid \beta;\mu^1,
\cdots,\mu^k)}
\newcommand{\Mbar}{\overline{\cM}}
\newcommand{\MgX}{\Mbar_{g,0}(X,\mu^+,\mu^-)}
\newcommand{\MX}{\Mbar^\bullet_\chi(X,\mu^+,\mu^-)}
\newcommand{\GX}{G^\bullet_\chi(X,\mu^+,\mu^-)}
\newcommand{\MgP}{\Mbar_{g,0}(\bP^1,\mu^+,\mu^-)}
\newcommand{\MP}{\Mbar^\bullet_\chi(\bP^1,\mu^+,\mu^-)}
\newcommand{\GP}{G^\bullet_\chi(\bP^1,\mu^+,\mu^-)}
\newcommand{\Mmn}{\Mbar^\bullet_{\chi^+}(\bP^1,\mu^+,\nu)}
\newcommand{\Mnm}{\Mbar^\bullet_{\chi^-}(\bP^1,\nu,\mu^-)}
\newcommand{\Mp}{\Mbar^\bullet_{\chi^+}(\bP^1,\mu^+,\nu^+)}
\newcommand{\Mm}{\Mbar^\bullet_{\chi^-}(\bP^1,\nu^-,\mu^-)}
\newcommand{\domain}{(C, \{ x_i^1\}_{i=1}^{l(\mu^1)},\ldots, \{ x_i^k\}_{i=1}^{l(\mu^k)} )}
\newcommand{\Cxy}{(C,\{x_i\}_{i=1}^{l(\mu^+)},\{y_j\}_{j=1}^{l(\mu^-)})}

\newcommand{\gmu}[1]{\mathcal{#1}_{g,\mu^+,\mu^-}}

\newcommand{\gm}[1]{{#1}_{g,\mu^+,\mu^-}}
\newcommand{\xm}[1]{{#1}^\bullet_{\chi,\mu^+,\mu^-}}
\newcommand{\xnx}[1]{{#1}_{\chi^+,\nu,\chi^-} }
\newcommand{\xmm}[1]{{#1}_{\chi,\mu^-,2l(\mu^-)}}
\newcommand{\mmx}[1]{{#1}_{2l(\mu^+),\mu^+,\chi}}
\newcommand{\xnxnx}[1]{{#1}_{\chi^+,\nu^+,\chi^0,\nu^-, \chi^-} }

\newcommand{\con}{ {2g-2 +l(\mu^+)+l(\mu^-) }}
\newcommand{\dis}{ {-\chi+l(\mu^+)+l(\mu^-) }}
\newcommand{\xmn}{ {-\chi^+ +l(\mu^+) + l(\nu) }}
\newcommand{\xnm}{ {-\chi^- +l(\nu)+ l(\mu^-)}}
\newcommand{\nmp}{ {-\chi^+ +l(\nu^+) + l(\mu^+) }}
\newcommand{\nmm}{ {-\chi^- +l(\nu^-) + l(\mu^-)}}

\newcommand{\two  }[1]{  ({#1}_1,{#1}_2)  }

\newcommand{\three}[1]{  ({#1}_1,{#1}_2,{#1}_3)  }

\newcommand{\five }[1]{  ({#1}_1,{#1}_2,{#1}_3,{#1}_4,{#1}_5)  }
\newcommand{\Five }[1]{ {[{#1}_1,{#1}_2,{#1}_3,{#1}_4,{#1}_5]} }
\newcommand{\tri}{(\chi^+,\nu,\chi^-)}
\newcommand{\lab}{ (\chi^+,\nu^+,\chi^0,\nu^-,\chi^-) }
\newcommand{\mm}{{\mu^+,\mu^-}}
\newcommand{\nn}{{\nu^+,\nu^-}}
\newcommand{\mn}{{\mu^+,\nu}}
\newcommand{\nm}{{\nu,\mu^-}}

\newcommand{\lam}{\lambda}
\newcommand{\si}{\sigma}

\newcommand{\pa}{\partial}
\newcommand{\bu}{\bullet}
\newcommand{\da}{D^\alpha }
\newcommand{\Da}{\Delta(D^\alpha)}
\newcommand{\dpm}{\Delta^\pm}

\newcommand{\TY}{f^*\left(\Omega_{Y[m^1,\ldots,m^k]}
(\sum_{\alpha=1}^k \log D^\alpha_{(m^\alpha)})\right)^\vee }
\newcommand{\zero}{H_{\mathrm{et}}^0(\mathbf{R}_l^{\alpha\bu})}
\newcommand{\one}{H_{\mathrm{et}}^1(\mathbf{R}_l^{\alpha\bu})}
\newcommand{\Ym}{Y[m^1,\ldots,m^k]}
\newcommand{\pim}{\pi[m^1,\ldots,m^k]}
\newcommand{\pp}{p^+_{\mu^+}p^-_{\mu^-}}
\newcommand{\amm}{|\Aut(\mu^+)||\Aut(\mu^-)|}

\DeclareMathOperator{\Aut}{Aut}
\DeclareMathOperator{\Ext}{Ext}

\DeclareMathOperator{\Br}{Br}
\DeclareMathOperator{\Sym}{Sym}

\newcommand{\ft}{\mathfrak t}

\newtheorem{theorem}{Theorem}[section]
\newtheorem{Theorem}{Theorem}
\newtheorem{proposition}[theorem]{Proposition}
\newtheorem{lemma}[theorem]{Lemma}
\newtheorem{rem}[theorem]{Remark}
\newtheorem{df}[theorem]{Definition}
\newtheorem{corollary}[theorem]{Corollary}
\newtheorem{ex}[theorem]{Example}

\theoremstyle{remark}

\theoremstyle{definition}

\definecolor{light-pink}{rgb}{1,.90,.90}
\definecolor{light-green}{rgb}{.95,1,.95}

\definecolor{yellow}{rgb}{1,1,0}
\definecolor{orange}{rgb}{1,.7,0}
\definecolor{red}{rgb}{1,0,0}
\definecolor{white}{rgb}{1,1,1}

\definecolor{A}{rgb}{.75,1,.75}

\begin{document}

\title{A Formula of Two-Partition Hodge Integrals}
\author{Chiu-Chu Melissa Liu}
\address{Department of Mathematics,
Harvard University, Cambridge, MA 02138, USA}
\email{ccliu@math.harvard.edu}
\author{Kefeng Liu}
\address{Center of Mathematical Sciences, Zhejiang University, Hangzhou, Zhejiang 310027, China;
Department of Mathematics,University of California at Los Angeles, Los Angeles, CA 90095-1555, USA\\
}
\email{liu@cms.zju.edu.cn, liu@math.ucla.edu}
\author{Jian Zhou}
\address{Department of Mathematical Sciences,
Tsinghua University, Beijing, 100084, China}
\email{jzhou@math.tsinghua.edu.cn}

\maketitle

\section{Introduction} \label{sec:Introduction}

Let $\Mbar_{g,n}$ denote the Deligne-Mumford moduli stack of stable curves of genus
$g$ with $n$ marked points. Let $\pi:\Mbar_{g,n+1}\to \Mbar_{g,n}$
be the universal curve, and let $\omega_\pi$ be the relative dualizing sheaf.
The Hodge bundle
$$
\bE=\pi_*\omega_\pi
$$
is a rank $g$ vector bundle over $\Mbar_{g,n}$ whose
fiber over $[(C,x_1,\ldots,x_n)]\in\Mbar_{g,n}$ is
$H^0(C,\omega_C)$.
Let $s_i:\Mbar_{g,n}\to\Mbar_{g,n+1}$ denote the section of $\pi$ which
corresponds to the $i$-th marked point, and let
$$
\bL_i=s_i^*\omega_\pi
$$
be the line bundle over $\Mbar_{g,n}$ whose fiber over
$[(C,x_1,\ldots,x_n)]\in\Mbar_{g,n}$ is the cotangent line $T^*_{x_i} C$
at the $i$-th marked point $x_i$.
A Hodge integral is an integral of the form
$$\int_{\Mbar_{g, n}} \psi_1^{j_1}
\cdots \psi_n^{j_n}\lam_1^{k_1} \cdots \lam_g^{k_g}$$
where
$\psi_i=c_1(\bL_i)$ is the first Chern class of $\bL_i$, and
$\lam_j=c_j(\bE)$ is the $j$-th Chern class of the Hodge bundle.

The study of Hodge integrals is an important part of the
intersection theory on $\Mbar_{g,n}$.  Hodge integrals also
naturally arise when one computes Gromov-Witten invariants by
localization techniques. For example, the following generating
series of Hodge integrals arises when one computes local invariants
of a toric Fano surface in a Calabi-Yau 3-fold by virtual
localization \cite{Zho2}:
\begin{equation} \label{eqn:Two}
\begin{split}
& G_{\mu^+, \mu^-}(\lambda; \tau)
=  - \frac{(\sqrt{-1}\lambda)^{l(\mu^+)+l(\mu^-)}}{z_{\mu^+} \cdot z_{\mu^-}}
\left[\tau(\tau+1)\right]^{l(\mu^+)+l(\mu^-)-1} \\
& \cdot \prod_{i=1}^{l(\mu^+)} \frac{\prod_{a=1}^{\mu^+_i-1}
\left( \mu^+_i\tau + a \right)}{\mu_i^+!}
\cdot \prod_{i=1}^{l(\mu^-)} \frac{\prod_{a=1}^{\mu^-_i-1}
\left( \mu_i^- \frac{1}{\tau} + a \right)}{\mu_i^-!} \\
&\cdot  \sum_{g \geq 0} \lambda^{2g-2} \int_{\Mbar_{g, l(\mu^+)+l(\mu^-)}}
\frac{\Lambda_{g}^{\vee}(1)\Lambda^{\vee}_{g}(\tau)\Lambda_{g}^{\vee}(-\tau - 1)}
{\prod_{i=1}^{l(\mu^+)} \frac{1}{\mu_i^+} \left(\frac{1}{\mu^+_i} - \psi_i\right)
\prod_{j=1}^{l(\mu^-)} \frac{\tau}{\mu_i^-}\left(\frac{\tau}{\mu^-_j} - \psi_{l(\mu^+)+j}\right)},
\end{split} \end{equation}
where $\lam, \tau$ are variables, $(\mu^+,\mu^-)\in \cP^2_+$, the
set of pairs of partitions which are not both empty, and
$$
\Lambda_g^\vee(u)=u^g-\lam_1 u^{g-1}+\cdots+(-1)^g\lam_g.
$$
We will call the Hodge integrals in $G_{\mu^+,\mu^-}(\lam;\tau)$ the
{\em two-partition Hodge integrals}.

The purpose of this paper is to prove the following formula conjectured in \cite{Zho3}:
\begin{equation} \label{eqn:Conj}
G^\bu(\lam;p^+,p^-;\tau)=R^\bu(\lam;p^+,p^-;\tau)
\end{equation}
where
$$
G^\bu(\lam;p^+,p^-;\tau)
= \exp \left(
\sum_{(\mu^+, \mu^-) \in \cP^2_+} G_{\mu^+,
  \mu^-}(\lam;\tau)p^+_{\mu^+}p^-_{\mu^-}\right)
$$

\begin{eqnarray*}
&&R^\bu(\lam;p^+,p^-;\tau)\\
&=& \sum_{|\mu^{\pm}|=|\nu^{\pm}|}
\frac{\chi_{\nu^+}(\mu^+)}{z_{\mu^+}} \frac{\chi_{\nu^-}(\mu^-)}{z_{\mu^-}}
e^{\sqrt{-1}(\kappa_{\nu^+} \tau  + \kappa_{\nu^-} \tau^{-1})\lambda/2}
\cW_{\nu^+, \nu^-}(e^{\sqrt{-1}\lam}) p^+_{\mu^+}p^-_{\mu^-},
\end{eqnarray*}
$p^\pm=(p^\pm_1,p^\pm_2,\ldots)$ are formal variables, and
$$
p^\pm_\mu=p^\pm_{\mu_1}\cdots p^\pm_{\mu_h}
$$
if $\mu=(\mu_1\geq\cdots\geq\mu_h>0)$. See Section \ref{sec:Conjecture} for
notation  in the definition of $R^\bu(\lam;p^+,p^-;\tau)$.

Formula (\ref{eqn:Conj}) is motivated by a formula of one-partition
Hodge integrals conjectured by M. Mari\~{n}o and C. Vafa in \cite{Mar-Vaf}
and proved by us in \cite{LLZ}. See \cite{Oko-Pan} for another approach
to the Mari\~{n}o-Vafa formula. The Mari\~{n}o-Vafa formula can be obtained
by setting $p^-=0$ in (\ref{eqn:Conj}). In a recent paper \cite{Dia-Flo}, 
D.E. Diaconescu and  B. Florea conjectured a relation between 
three-partition Hodge integrals and the topological vertex
\cite{AKMV}. A mathematical theory of the topological vertex will be 
developed in \cite{LLLZ}.

The generating function $R^\bu(\lam;p^+,p^-;\tau)$
is a combinatorial expression involving the representation theory of
Kac-Moody Lie algebras. It is also related to the HOMFLY polynomial of the
Hopf link and the Chern-Simon theory \cite{Wit1, Mor-Luk}. 
In \cite{Zho4}, the third author used (\ref{eqn:Conj}) and a combinatorial
trick called the chemistry of $\bZ_k$-colored labelled graphs to
prove a formula conjectured by A. Iqbal in \cite{Iqb} which expresses
the generating function of Gromov-Witten invariants in all genera of
local toric Calabi-Yau threefolds in terms of $\cW_{\mu,\nu}$. 
See \cite{Iqb, AKMV, Dia-Flo} for surveys of works on this subject. 

Our strategy to prove (\ref{eqn:Conj}) is based on the following
cut-and-join equation of $R^\bu$ observed in \cite{Zho3}:
\begin{equation}\label{eqn:Rcj}
\frac{\pa}{\pa\tau}R^\bu=\frac{\sqrt{-1}\lam}{2}(C^+ + J^+)R^\bu 
-\frac{\sqrt{-1}\lam}{2\tau^2}(C^- +
J^-)R^\bu
\end{equation}
where
$$
C^\pm=\sum_{i,j}(i+j)p_i^\pm p_j^\pm\frac{\pa}{\pa p^\pm_{i+j}},\ \ \ \ \
J^\pm=\sum_{i,j}ij p_{i+j}^\pm\frac{\pa^2}{\pa p_i^\pm\pa p_j^\pm}.
$$
Equation (\ref{eqn:Rcj}) can be derived by the method in \cite{Zho1, LLZ}. 
In \cite{Zho3}, the third author proved that
\begin{Theorem}[initial values]\label{thm:initial}
\begin{equation}\label{eqn:initial}
G^\bu(\lam;p^+,p^-;-1)=R^\bu(\lam;p^+,p^-;-1).
\end{equation}
\end{Theorem}
So (\ref{eqn:Conj}) follows from the main theorem in this paper:
\begin{Theorem}[cut-and-join equation of $G^\bu$]\label{thm:Gcj}
\begin{equation}\label{eqn:Gcj}
\frac{\pa}{\pa\tau}G^\bu=
\frac{\sqrt{-1}\lam}{2}(C^+ + J^+)G^\bu 
-\frac{\sqrt{-1}\lam}{2\tau^2}(C^- + J^-)G^\bu
\end{equation}
\end{Theorem}

Both \cite[Theorem 2]{LLZ} (cut-and-join equation of one-partition
Hodge integrals) and Theorem \ref{thm:Gcj} are proved by
localization method. We compute certain relative Gromov-Witten
invariants by virtual localization, and get an expression in terms
of one-partition or two-partition Hodge integrals and certain
integrals of target $\psi$ classes. In \cite{LLZ}, we used
functorial localization to push forward calculations to projective
spaces, where the equivariant cohomology is completely understood,
and derived \cite[Theorem 2]{LLZ} without using much information
about integrals of target $\psi$ classes. In this paper, we relate
integrals of target $\psi$ classes to double Hurwitz numbers, and
use properties of double Hurwitz numbers to prove Theorem
\ref{thm:Gcj}. More precisely, for each $(\mm)\in \cP^2_+$, we
will define a generating function
$$
K^\bu_\mm(\lam)
$$
of certain relative Gromov-Witten invariants of $\bP^1\times\bP^1$
blowup at a point, and use localization method to derive the
following expression:
\begin{equation}\label{eqn:KG}
K^\bu_\mm(\lam)=\sum_{|\nu^\pm|=|\mu^\pm|}
\Phi^\bu_{\mu^+,\nu^+}(-\sqrt{-1}\tau\lam)z_{\nu^+}G^\bu_\nn(\lam;\tau)
z_{\nu^-}\Phi^\bu_{\nu^-,\mu^-}(\frac{-\sqrt{-1}}{\tau}\lam)
\end{equation}
In (\ref{eqn:KG}), $\Phi_{\mu,\nu}^\bu(\lam)$ is a generating function of double Hurwitz numbers,
and $z_\mu$ is defined in Section \ref{sec:partition}. It turns out that
(\ref{eqn:KG}) is equivalent to the following equation:
\begin{equation}\label{eqn:GK}
G^\bu_\mm(\lam;\tau)=\sum_{|\nu^\pm|=|\mu^\pm|}
\Phi^\bu_{\mu^+,\nu^+}(\sqrt{-1}\tau\lam)z_{\nu^+}K^\bu_\nn(\lam)
z_{\nu^-}\Phi^\bu_{\nu^-,\mu^-}(\frac{\sqrt{-1}}{\tau}\lam)
\end{equation}
So Theorem \ref{thm:Gcj} (cut-and-join equation of $G^\bu$)
follows from the cut-and-join equations of double Hurwitz numbers.
As a consequence, one can compute $K_{\mu^+, \mu^-}(\lambda)$ in
terms of $\cW_{\nu^+, \nu^-}$ (Corollary \ref{zeroframing}). 
We will give three derivations of the cut-and-join equations
of double Hurwitz numbers: by combinatorics (Section \ref{sec:Hcj}), by
gluing formula (Section \ref{sec:glue}), and by localization (Section
\ref{sec:ELSV}).

The rest of the paper is arranged as follows. In Section
\ref{sec:Conjecture}, we give the precise statement of
(\ref{eqn:Conj}), and recall the proof of Theorem
\ref{thm:initial} (initial values). In Section \ref{sec:Hurwitz}
we give a combinatorial study of double Hurwitz numbers, and
derive Theorem \ref{thm:Gcj} (the cut-and-join equation of
$G^\bu$) from (\ref{eqn:KG}) and some identities of double Hurwitz
numbers. In Section \ref{sec:LiGV}, we review J. Li's works
\cite{Li1, Li2} on moduli spaces of relative stable morphisms, and
virtual localization on such moduli spaces \cite{Gra-Pan, Gra-Vak2}. 
In Section \ref{sec:LocalizeP}, we give a geometric study of double Hurwitz
numbers. In Section \ref{sec:Objects}, we introduce the geometric
objects involved  in the proof of (\ref{eqn:KG}). In Section
\ref{sec:LocalizeX}, we prove (\ref{eqn:KG}) by arranging the
localization contribution in a neat way.

\bigskip
\noindent
{\bf Acknowledgments.}
We wish to thank Jun Li for explaining his works \cite{Li1, Li2} and
Ravi Vakil for explaining relative virtual localization \cite{Gra-Vak2}.
The research in this work was started during the visit of the first and the third authors
to the Center of Mathematical Sciences, Zhejiang University in July and August of 2003.
The hospitality of the Center is greatly appreciated.
The second author is supported by an NSF grant.
The third author is partially supported by research
grants from NSFC and Tsinghua University.

\pagebreak

\section{The Conjecture}
\label{sec:Conjecture}

\subsection{Partitions}\label{sec:partition}
We recall some notation of partitions. Given a partition
$$
\mu=(\mu_1\geq \mu_2\geq \cdots \geq \mu_h>0),
$$
write $l(\mu)=h$, and $|\mu|=\mu_1+\cdots+\mu_h$. Define
$$
\kappa_\mu=\sum_{i=1}^{l(\mu)}\mu_i(\mu_i-2i+1).
$$
For each positive integer $j$, define
$$
m_j(\mu)=|\{i:\mu_i=j\}|.
$$
Then
$$
|\Aut(\mu)|=\prod_j (m_j(\mu))!.
$$
Define
$$
z_\mu=\mu_1\cdots\mu_{l(\mu)}|\Aut(\mu)|=\prod_j\left(m_j(\mu)! j^{m_j(\mu)}\right).
$$

Let $\cP$ denote the set of partitions.
We allow the empty partition and take
$$l(\emptyset)=|\emptyset|=\kappa_{\emptyset} = 0.$$
Let
$$
\cP^2_+=\cP^2-\{(\emptyset,\emptyset)\}.
$$

\subsection{Generating functions of two-partition Hodge integrals}
\label{sec:hodge}
For $(\mu^+, \mu^-) \in \cP^2_+$, define
\begin{eqnarray*}
&& \gm{G}(\alpha,\beta)\\
&=&\frac{ -\sqrt{-1}^{l(\mu^+)+l(\mu^-)} }{|\Aut(\mu^+)||\Aut(\mu^-)|}
 \prod_{i=1}^{l(\mu^+)}
\frac{\prod_{a=1}^{\mu^+_i-1}(\mu^+_i\beta + a\alpha) }
     {(\mu^+_i-1)!\alpha^{\mu^+_i-1} }
\prod_{j=1}^{l(\mu^-)}
\frac{\prod_{a=1}^{\mu^-_j-1}(\mu^-_j\alpha + a\beta)}
     {(\mu^-_j-1)!\beta^{\mu^-_j-1} }\\
&&\cdot \int_{\Mbar_{g,l(\mu^+)+l(\mu^-)}}
\frac{\Lambda^\vee(\alpha)\Lambda^\vee(\beta)
\Lambda^\vee(-\alpha-\beta)(\alpha\beta(\alpha+\beta))^{l(\mu^+)+l(\mu^-)-1}}
{\prod_{i=1}^{l(\mu^+)}(\alpha(\alpha-\mu^+_i\psi_i))
 \prod_{j=1}^{l(\mu^-)}(\beta (\beta- \mu^-_j\psi_{l(\mu^+)+j})}.
\end{eqnarray*}
We have the following special cases which have been studied in \cite{LLZ}:
\begin{eqnarray*}
G_{g,\mu^+,\emptyset}(\alpha,\beta)
&=&\frac{ -\sqrt{-1}^{l(\mu^+)} }{|\Aut(\mu^+)|}
 \prod_{i=1}^{l(\mu^+)}
\frac{\prod_{a=1}^{\mu^+_i-1}(\mu^+_i\beta + a\alpha) }
     {(\mu^+_i-1)!\alpha^{\mu^+_i-1} }\\
&&\cdot \int_{\Mbar_{g,l(\mu^+)}}
\frac{\Lambda^\vee(\alpha)\Lambda^\vee(\beta)
\Lambda^\vee(-\alpha-\beta)(\alpha\beta(\alpha+\beta))^{l(\mu^+)-1}}
{\prod_{i=1}^{l(\mu^+)}(\alpha(\alpha-\mu^+_i\psi_i))}\\
G_{g,\emptyset,\mu^-}(\alpha,\beta)
&=&\frac{ -\sqrt{-1}^{l(\mu^-)} }{|\Aut(\mu^-)|}
\prod_{j=1}^{l(\mu^-)}
\frac{\prod_{a=1}^{\mu^-_j-1}(\mu^-_j\alpha + a\beta)}
     {(\mu^-_j-1)!\beta^{\mu^-_j-1} }\\
&&\cdot \int_{\Mbar_{g,l(\mu^-)}}
\frac{\Lambda^\vee(\alpha)\Lambda^\vee(\beta)
\Lambda^\vee(-\alpha-\beta)(\alpha\beta(\alpha+\beta))^{l(\mu^-)-1}}
{\prod_{j=1}^{l(\mu^-)}(\beta (\beta- \mu^-_j\psi_j))}
\end{eqnarray*}
By a standard degree argument,
one sees that $\gm{G}(\alpha,\beta)$ is homogeneous of degree $0$, so
$$
\gm{G}(\alpha,\beta)=\gm{G}(1,\frac{\beta}{\alpha}).
$$
Let
$$
\gm{G}(\tau)=\gm{G}(1,\tau).
$$

Introduce variables $\lam$, $p^+=(p^+_1,p^+_2,\ldots)$,
$p^-=(p^-_1,p^-_2,\ldots)$. Given a partition $\mu$,
define
$$
p^\pm_\mu=p^\pm_1\cdots p^\pm_{l(\mu)}.
$$
In particular, $p^\pm_\emptyset=1$. Define
\begin{eqnarray*}
G_\mm(\lam;\tau)&=&\sum_{g=0}^\infty \lam^\con\gm{G}(\tau)\\
G(\lam;p^+,p^-;\tau)&=&\sum_{(\mm)\in\cP^2_+}
G_\mm(\lam;\tau)\pp\\
G^\bu(\lam;p^+,p^-;\tau)&=&\exp(G(\lam;p^+,p^-;\tau))\\
&=& \sum_{(\mm)\in\cP^2} G^\bu_\mm(\lam;\tau)\pp\\
G^\bu_\mm(\lam;\tau)&=& \sum_{\chi\in2\bZ,\chi\leq 2(l(\mu^+)+l(\mu^-))}
\lam^{\dis}\xm{G}(\tau)
\end{eqnarray*}

\subsection{Generating functions of representations of symmetric groups}
Let
$$
q=e^{\sqrt{-1}\lambda},\ \ \ [m]=q^{m/2}-q^{-m/2}.
$$
Define
\begin{eqnarray}
&& \cW_{\mu, \nu}(q) = q^{|\nu|/2} \cW_{\mu}(q) \cdot s_{\nu}(\cE_{\mu}(t)),
\end{eqnarray}
where
\begin{eqnarray}
&& \cW_{\mu}(q) = q^{\kappa_{\mu}/4}\prod_{1 \leq i < j \leq l(\mu)}
\frac{[\mu_i - \mu_j + j - i]}{[j-i]}
\prod_{i=1}^{l(\mu)} \prod_{v=1}^{\mu_i} \frac{1}{[v-i+l(\mu)]}, \\
&& \cE_{\mu}(t) = \prod_{j=1}^{l(\mu)} \frac{1+q^{\mu_j-j}t}{1+q^{-j}t}
\cdot \left(1 + \sum_{n=1}^{\infty}
\frac{t^n}{\prod_{i=1}^n (q^i-1)}\right).
\end{eqnarray}
In the special case of $(\mu^+, \mu^-) = (\emptyset, \emptyset)$,
we have
$$\cW_{\emptyset, \emptyset} = 1.$$

Define
\begin{eqnarray*}
R^\bu(\lam;p^+,p^-;\tau)
& = &\sum_{|\nu^\pm|=|\mu^\pm| \geq 0}
\frac{\chi_{\nu^+}(C(\mu^+))}{z_{\mu^+}}
\frac{\chi_{\nu^-}(C(\mu^-))}{z_{\mu^-}} \\
&& \cdot e^{\sqrt{-1}(\kappa_{\nu^+}\tau+ \kappa_{\nu^-}\tau^{-1})\lambda/2}
\cW_\nn(e^{\sqrt{-1}\lam}) \pp.
\end{eqnarray*}

\subsection{The conjecture and the strategy}
The main purpose of this paper is to prove the following formula  conjectured
by the third author in \cite{Zho3}:
\begin{Theorem}
We have the following formula of two-partition Hodge integrals:

\bigskip
\noindent
\textup{(\ref{eqn:Conj})}\makebox[1in]{ }
$G^\bu(\lam;p^+, p^-;\tau) = R^\bu(\lam;p^+, p^-;\tau).$
\end{Theorem}

The method in \cite{Zho1, LLZ} shows that $R^\bu$ satisfies the cut-and-join
equation (\ref{eqn:Rcj}). In \cite{Zho3}, the third author proved Theorem
\ref{thm:initial} (initial values).
So (\ref{eqn:Conj}) follows from Theorem \ref{thm:Gcj} (cut-and-join equation
of $G^\bu$). We will recall the proof of  Theorem \ref{thm:initial}
in Section \ref{sec:initial}.

\subsection{Initial values} \label{sec:initial}
For completeness, we now recall the proof of Theorem \ref{thm:initial}, which
says

\bigskip
\noindent
(\ref{eqn:initial})\makebox[1in]{ }
$G^\bu(\lam;p^+, p^-;-1) = R^{\bu}(\lam;p^+, p^-;-1).$
\bigskip

We need the skew Schur functions \cite{Mac}.
Recall the Schur functions are related to the Newton functions by:
$$
s_\mu(x)=\sum_{|\nu|=|\mu|}\frac{\chi_\mu(\nu)}{z_\nu} p_\nu(x),
$$
where $x=(x_1,x_2,\ldots)$ are formal variables such that
$$
p_i(x)=x_1^i + x_2^i+\cdots.
$$
There are integers $c_{\mu\nu}^{\eta}$ such that
$$s_{\mu}s_{\nu} = \sum_{\eta} c^{\eta}_{\mu\nu} s_{\eta}.$$
The skew Schur functions are defined by:
$$s_{\eta/\mu} = \sum_{\nu} c^{\eta}_{\mu\nu}s_{\nu}.$$
Note that $p^\pm= p(x^\pm)$.

\subsubsection{The left-hand-side}
When $l(\mu^+)+l(\mu^-) > 2$,
$$
G_\mm(\lam; -1) = 0;
$$
when $l(\mu^+) = 1$ and $l(\mu^-) = 0$,
\begin{eqnarray*}
&& G_\mm(\lam; -1) \\
& = &  - \sqrt{-1} \lam^{-1}  \sum_{g \geq 0} \lambda^{2g}
\int_{\Mbar_{g, 1}} \frac{\lam_g} {\frac{1}{\mu_1^+}
\left(\frac{1}{\mu^+_1} - \psi_1\right)}
\frac{\prod_{a=1}^{\mu^+_1-1} \left(-\mu_1^++ a \right)}{\mu_1^+ \cdot \mu_1^+!} \\
& = &  (-1)^{\mu^+_1} \sqrt{-1} \cdot \frac{1}{2\mu_1^+\sin (\mu_1^+\lambda/2)}
= \frac{(-1)^{\mu_1^+-1}}{q^{\mu_1^+/2} - q^{-\mu_1^+/2}} \cdot \frac{p_{\mu_1^+}}{\mu_1^+};
\end{eqnarray*}
the case of $l(\mu^+) = 0$ and $l(\mu^-) = 1$ is similar;
when $l(\mu^+) = l(\mu^-) = 1$,
\begin{eqnarray*}
G_\mm(\lam; -1)
& = &  \lim_{\tau \to -1}  \sum_{g \geq 0} \lambda^{2g} \int_{\Mbar_{g, 2}}
\frac{\Lambda_{g}^{\vee}(1)\Lambda^{\vee}_{g}(\tau)\Lambda_{g}^{\vee}(-1 - \tau)}
{\frac{1}{\mu_1^+} \left(\frac{1}{\mu^+_1} - \psi_1\right)
\cdot \frac{\tau}{\mu^-_1} \left(\frac{\tau}{\mu^-_1} - \psi_{2}\right)} \\
&& \cdot \tau(1+\tau) \cdot
\frac{\prod_{a=1}^{\mu^+_1-1} \left(\mu_1^+ \tau + a \right)}{\mu_1^+ \cdot \mu_1^+!}
\cdot \frac{\prod_{a=1}^{\mu^-_1-1} \left(\frac{\mu_1^-}{\tau} + a \right)}{\mu_1^- \cdot \mu_1^-!}.
\end{eqnarray*}
One needs to consider the $g=0$ term and the $g>0$ terms separately.
In the second case, the limit is zero while in first case,
by our convention:
$$\int_{\Mbar_{0, 2}}
\frac{\Lambda_{0}^{\vee}(1)\Lambda^{\vee}_{0}(\tau)\Lambda_{0}^{\vee}(-1 - \tau)}
{\frac{1}{\mu_1^+} \left(\frac{1}{\mu^+_1} - \psi_1\right)
\cdot \frac{\tau}{\mu^-_1} \left(\frac{\tau}{\mu^-_1} - \psi_{2}\right)}
= \frac{(\mu_1^+)^2(\frac{\mu_1^-}{\tau})^2}{\mu_1^+ + \frac{\mu_1^-}{\tau}},$$
hence when $\mu_1^+ \neq \mu_1^-$,
the limit is zero,
when $\mu_1^+ = \mu_1^-$,
the limit is:
\begin{eqnarray*}
&& \lim_{\tau \to -1}
\frac{(\mu_1^+)^2(\frac{\mu_1^-}{\tau})^2}{\mu_1^+ + \frac{\mu_1^-}{\tau}}
\cdot
\tau(1+\tau) \cdot
\frac{\prod_{a=1}^{\mu^+_1-1} \left(\mu_1^+ \tau + a \right)}{\mu_1^+ \cdot \mu_1^+!}
\cdot \frac{\prod_{a=1}^{\mu^-_1-1} \left(\frac{\mu_1^-}{\tau} + a \right)}{\mu_1^-\cdot\mu_1^-!}
= \frac{1}{\mu_1^+}.
\end{eqnarray*}

Recall that $p^\pm=p(x^\pm)$. With this notation, the initial value is:
\begin{eqnarray*}
&& G^\bu(\lam;p(x^+),p(x^-);-1) \\
& = & \exp \left(\sum_{n \geq 1} \frac{(-1)^{n-1}}{q^{n/2} - q^{-n/2}}\frac{p_n(x^+)}{n}
+ \sum_{n \geq 1} \frac{(-1)^{n-1}}{q^{n/2} - q^{-n/2}}\frac{p_n(x^-)}{n}
+ \sum_{n \geq 1} \frac{p_n(x^+)p_n(x^-)}{n}\right) \\
& = & \prod_{i,j=1}^{\infty} \frac{1}{(1 + q^{i-1/2}x^+_j)(1+q^{i-1/2}x^-_j)}
\prod_{j,k} \frac{1}{1-x^+_j x^-_k} \\
& = & \sum_{\rho^+} s_{\rho^+}(-q^{1/2}, -q^{3/2}, \dots)s_{\rho^+}(x^+)
\cdot \sum_{\rho} s_{\rho}(x^+)s_{\rho}(x^-) \\
&& \cdot \sum_{\nu^-} s_{\rho^-}(-q^{1/2}, -q^{3/2}, \dots)s_{\rho^-}(x^-) \\
& = & \sum_{\nu^{\pm}, \rho, \rho^{\pm}}
 s_{\rho^+}(-q^{1/2}, -q^{3/2}, \dots)c_{\rho^+\rho}^{\nu^+}s_{\nu^+}(x^+)
\cdot c_{\rho^-\rho}^{\nu^-} s_{\nu^-}(x^-) s_{\rho^-}(-q^{1/2}, -q^{3/2}, \dots) \\
& = & \sum_{\rho, \nu^{\pm}}
s_{\nu^+/\rho}(-q^{1/2}, -q^{3/2}, \dots)s_{\nu^-/\rho}(-q^{1/2}, -q^{3/2}, \dots)
\cdot s_{\nu^+}(x^+) s_{\nu^-}(x^-).
\end{eqnarray*}

\subsubsection{The right-hand side}

The following identity is proved in \cite{Zho3}:
\begin{eqnarray} \label{eqn:Key}
&& \cW_{\mu, \nu}(q)
= (-1)^{|\mu|+|\nu|}
q^{\frac{\kappa_{\mu}+\kappa_{\nu}+|\mu|+|\nu|}{2}}
\sum_{\rho} q^{-|\rho|} s_{\mu/\rho}(1, q, \dots)s_{\nu/\rho}(1, q, \dots).
\end{eqnarray}
From this one gets:
\begin{eqnarray*}
&& R^{\bullet}(\lambda; p(x^+), p(x^-);-1) \\
& = & \sum_{|\mu^{\pm}|=|\nu^{\pm}|}
\frac{\chi_{\nu^+}(\mu^+)}{z_{\mu^+}}
\frac{\chi_{\nu^-}(\mu^-)}{z_{\mu^-}}
e^{-\sqrt{-1}(\kappa_{\nu^+}+ \kappa_{\nu^-})\lambda/2} \cW_{\nu^+, \nu^-}(q)
p^+_{\mu^+}p^-_{\mu^-} \\
& = & \sum_{\nu^{\pm}} s_{\nu^+}(x^+) q^{-\kappa_{\nu^+}/2}
\cW_{\nu^+, \nu^-}(q) q^{-\kappa_{\nu^-}/2}s_{\nu^-}(x^-) \\
& = & \sum_{\nu^{\pm}} s_{\nu^+}(x^+)
s_{\nu^-}(x^-)
(-1)^{|\nu^+|+|\nu^-|}q^{(|\nu^+|+|\nu^-|)/2} \\
&& \cdot \sum_{\rho} q^{-|\rho|} s_{\nu^+/\rho}(1, q, \dots)s_{\nu^-/\rho}(1, q, \dots) \\
& = & \sum_{\nu^{\pm}} s_{\nu^+}(x^+)
s_{\nu^-}(x^-)
 \sum_{\rho} s_{\nu^+/\rho}(-q^{1/2}, -q^{3/2}, \dots)s_{\nu^-/\rho}(-q^{1/2}, -q^{3/2}, \dots).
\end{eqnarray*}
The proof of Theorem \ref{thm:initial} is complete.

\section{Double Hurwitz numbers and the cut-and-join equation of $G^\bu$}
\label{sec:Hurwitz}

In this section, we first derive some identities of double Hurwitz
numbers, such as sum formula and cut-and-join equations, which,
together with initial values, characterize the double Hurwitz
numbers. Then we combine these identities with (\ref{eqn:KG}) to
obtain Theorem \ref{thm:Gcj} (cut-and-join equation of $G^\bu$).

\subsection{Double Hurwitz numbers}

Let $X$ be a Riemann surface of genus $h$.
Given $n$ partitions $\eta^1, \dots, \eta^n$ of $d$,
denote by $H^X_d(\eta^1, \dots, \eta^n)^{\bullet}$ and $H^X_d(\eta^1, \dots, \eta^n)^{\circ}$
the weighted counts of possibly disconnected and connected Hurwitz covers of type $(\eta^1, \dots, \eta^n)$
respectively.
We will use the following formula for Hurwitz numbers (see e.g. \cite{Dij}):
\begin{eqnarray} \label{eqn:Burnside}
&& H^X_d(\eta^1, \dots, \eta^n)^{\bullet}
= \sum_{|\rho|= d} \left(\frac{\dim R_\rho}{d!}\right)^{2 - 2h}
\prod_{i=1}^n |C_{\eta^i}|\frac{\chi_{\rho}(C_{\eta^i})}{\dim R_{\rho}}.
\end{eqnarray}
It is sometimes referred to as the {\em Burnside formula}.

Suppose $C \to \bP^1$ is a genus $g$ cover
which has ramification type $\mu^+, \mu^-$ at two points $p_0$ and $p_1$ respectively,
and ramification type $(2)$ at $r$ other points.
By Riemann-Hurwitz formula,
\begin{eqnarray} \label{eqn:r}
r = 2g - 2 + l(\mu^+) + l(\mu^-).
\end{eqnarray}
Denote
\begin{eqnarray*}
&& H^\circ_g(\mm) = H^{\bP^1}_d(\mu^+, \mu^-, \eta^1, \dots, \eta^{r})^{\circ}, \\
&& H^\bu_g(\mm) = H^{\bP^1}_d(\mu^+, \mu^-, \eta^1, \dots, \eta^{r})^{\bullet},
\end{eqnarray*}
for $\eta^1 = \cdots = \eta^r = (2)$.
We have by (\ref{eqn:Burnside}):
\begin{eqnarray} \label{eqn:Simple}
H_g^{\bullet}(\mu^+, \mu^-)
& = & \sum_{|\nu| = d} f_\nu(2)^{r}
\frac{\chi_\nu(C_{\mu^+})}{z_{\mu^+}}\frac{\chi_{\nu}(C_{\mu^-})}{z_{\mu^+}},
\end{eqnarray}
where $r$ is given by (\ref{eqn:r}), and
$$f_{\nu}(2) = |C_{(2)}|\frac{\chi_{\nu}(C_{(2)})}{\dim R_{\nu}}.$$
Define
\begin{eqnarray*}
&& \Phi_{\mu^+, \mu^-}^{\circ}(\lambda)
= \sum_{g \geq 0} H_g^{\circ}({\mu^+, \mu^-})
\frac{\lambda^{2g-2+l(\mu^+)+l(\mu^-)}}{(2g-2 + l(\mu^+) + l(\mu^-))!}, \\
&& \Phi^{\bullet}_{\mu^+, \mu^-}(\lambda)
= \sum_{g \geq 0} H^{\bullet}_g({\mu^+, \mu^-})
\frac{\lambda^{2g-2+l(\mu^+)+l(\mu^-)}}{(2g-2 + l(\mu^+) + l(\mu^-))!}, \\
&& \Phi^{\circ}(\lambda; p^+, p^-)
= \sum_{\mu^+, \mu^-} \Phi^{\circ}_{\mu^+, \mu^-}(\lambda) p^+_{\mu^+}p_{\mu^-}^-, \\
&& \Phi^{\bullet}(\lambda; p^+, p^-)
= 1+\sum_{\mu^+, \mu^-} \Phi_{\mu^+, \mu^-}^{\bullet} (\lambda) p_{\mu^+}^+p_{\mu^-}^-.
\end{eqnarray*}
The usual relationship between connected and disconnected Hurwitz numbers is:
\begin{eqnarray} \label{eqn:Conn}
&&\Phi^{\circ}(\lambda; p^+, p^-) = \log \Phi^{\bullet}(\lambda; p^+, p^-).
\end{eqnarray}

By (\ref{eqn:Simple}) one easily gets:
\begin{eqnarray}
&& \Phi^{\bullet}(\lambda; p^+, p^-)
= 1+ \sum_{d \geq 1} \sum_{|\mu^\pm|= d} \sum_{|\nu|=d}
\frac{\chi_{\nu}(C_{\mu^+})}{z_{\mu^+}} \frac{\chi_{\nu}(C_{\mu^-})}{z_{\mu^-}}
e^{f_{\nu}(2)\lambda}
p_{\mu^+}^+p_{\mu^-}^-.
\end{eqnarray}
Equivalently,
\begin{eqnarray} \label{eqn:BUR1}
&& \Phi^{\bullet}_{\mu^+, \mu^-}(\lambda)
= \sum_{|\nu|=d}
\frac{\chi_{\nu}(C_{\mu^+})}{z_{\mu^+}} \frac{\chi_{\nu}(C_{\mu^-})}{z_{\mu^-}}
e^{f_{\nu}(2)\lambda}.
\end{eqnarray}
We also have
\begin{eqnarray} \label{eqn:BUR2}
&& \Phi^{\bullet}_{\mu^+, \mu^-}(\lambda)
= \sum_{|\nu|=d} (-1)^{l(\mu^+)+l(\mu^-)}
\frac{\chi_{\nu}(C_{\mu^+})}{z_{\mu^+}} \frac{\chi_{\nu}(C_{\mu^-})}{z_{\mu^-}}
e^{-f_{\nu}(2)\lambda}.
\end{eqnarray}

\subsection{Sum formula and initial values}\label{sec:sum}
\begin{proposition}
We have
\begin{eqnarray}
&&  \Phi^{\bullet}_{\mu^1, \mu^3}(\lambda_1+\lambda_2)
=  \sum_{\mu^2} \Phi^{\bullet}_{\mu^1, \mu^2}(\lambda_1)
\cdot z_{\mu^2} \cdot \Phi^{\bullet}_{\mu^2, \mu^3}(\lambda_2), \label{eqn:Comp}  \\
&& \Phi^{\bullet}_{\mu^1, \mu^3}(0) = \frac{1}{z_{\mu^1}}\delta_{\mu^1, \mu^3}. \label{eqn:Init}
\end{eqnarray}
\end{proposition}

\begin{proof}
By the orthogonality relation for characters of $S_d$:
\begin{eqnarray} \label{eqn:Orth1}
&& \sum_{\mu} \frac{\chi_{\nu^1}(C_{\mu})\chi_{\nu^2}(C_{\mu})}{z_{\mu}}
= \delta_{\nu^1, \nu^2}
\end{eqnarray}
we have
\begin{eqnarray*}
&& \sum_{\mu^2} \Phi^{\bullet}_{\mu^1, \mu^2}(\lambda_1)
\cdot z_{\mu^2} \cdot \Phi^{\bullet}_{\mu^2, \mu^3}(\lambda_2) \\
& = & \sum_{\mu^2} \sum_{\nu^1}
\frac{\chi_{\nu^1}(C_{\mu^1})}{z_{\mu^1}} \frac{\chi_{\nu^1}(C_{\mu^2})}{z_{\mu^2}}
e^{f_{\nu^1}(2)\lambda_1} \cdot z_{\mu^2} \cdot
\sum_{\nu^2}
\frac{\chi_{\nu^2}(C_{\mu^2})}{z_{\mu^2}} \frac{\chi_{\nu^2}(C_{\mu^3})}{z_{\mu^3}}
e^{f_{\nu^2}(2)\lambda_2} \\
& = & \sum_{\nu^1} \sum_{\nu^2}
\frac{\chi_{\nu^1}(C_{\mu^1})}{z_{\mu^1}} \frac{\chi_{\nu^2}(C_{\mu^3})}{z_{\mu^3}}
e^{f_{\nu^1}(2)\lambda_1+f_{\nu^2}(2)\lambda_2} \cdot
\sum_{\mu^2}
\frac{\chi_{\nu^1}(C_{\mu^2})\chi_{\nu^2}(C_{\mu^2})}{z_{\mu^2}} \\
& = &  \sum_{\nu^1} \sum_{\nu^2}
\frac{\chi_{\nu^1}(C_{\mu^1})}{z_{\mu^1}} \frac{\chi_{\nu^2}(C_{\mu^3})}{z_{\mu^3}}
e^{f_{\nu^1}(2)\lambda_1+f_{\nu^2}(2)\lambda_2} \delta_{\nu^1, \nu^2} \\
& = & \sum_{\nu}
\frac{\chi_{\nu}(C_{\mu^1})}{z_{\mu^1}} \frac{\chi_{\nu}(C_{\mu^3})}{z_{\mu^3}}
e^{f_{\nu}(2)(\lambda_1+\lambda_2)} \\
& = & \Phi^{\bullet}_{\mu^1, \mu^3}(\lambda_1+\lambda_2).
\end{eqnarray*}
Similarly,
by the orthogonality relation:
\begin{eqnarray} \label{eqn:Orth2}
&& \sum_{|\nu|=d} \chi_{\nu}(C_{\mu^1}) \cdot \chi_{\nu}(C_{\mu^2})
= z_{\mu^1}\delta_{\mu^1, \mu^2}.
\end{eqnarray}
we have
\begin{eqnarray*}
\Phi_{\mu^1, \mu^2}^{\bullet}(0)
& = & \sum_{|\nu|=d}  \frac{\chi_{\nu}(C_{\mu^1})}{z_{\mu^1}}
\cdot \frac{\chi_{\nu}(C_{\mu^2})}{z_{\mu^2}}
= \frac{1}{z_{\mu^1}}\delta_{\mu^1, \mu^2}.
\end{eqnarray*}
\end{proof}

Equation (\ref{eqn:Comp}) is a sum formula for double Hurwitz numbers, and
Equation (\ref{eqn:Init}) gives the initial values for double Hurwitz numbers.

\begin{corollary}
Denote by $\Phi^{\bullet}(\lambda)_d$ the matrix
$(\Phi^{\bullet}_{\mu, \nu}(\lambda))_{|\mu|=|\nu|=d}$.
Then $\Phi^{\bullet}(\lambda)_d$ is invertible,
and
\begin{eqnarray} \label{eqn:PhiInverse}
Z^{-1}_d \Phi^{\bullet}(-\lambda)_d^{-1} = \Phi^{\bullet}(\lambda)_d Z_d.
\end{eqnarray}
where $Z_d = (z_{\mu}\delta_{\mu, \nu})_{|\mu|=|\nu|=d}$.
\end{corollary}

\begin{proof}
In (\ref{eqn:Comp}) we take $\lambda_1 = \lambda$ and $\lambda_2 = - \lambda$,
then by (\ref{eqn:Init}) we have
$$Z^{-1}_d = \Phi^\bu(0)_d = \Phi^\bu(\lam)_d Z_d \Phi^\bu(-\lam)_d.$$
Taking determinant on both sides one sees that
$\Phi^\bu(\lam)_d$ is invertible,
and (\ref{eqn:PhiInverse}) is a straightforward consequence.
\end{proof}

\subsection{Cut-and-join equation for double Hurwitz numbers}
\label{sec:Hcj}
Recall for any partition $\nu$ of $d$, one has
\begin{eqnarray*}
&& f_{\nu}(2)  \cdot \sum_{\mu}  \frac{\chi_{\nu}(C_\mu)}{z_{\mu}}p_{\mu}
= \frac{1}{2} \sum_{i,j}
\left((i+j)p_ip_j \frac{\partial}{\partial p_{i+j}}
+ ijp_{i+j}\frac{\partial}{\partial p_i}\frac{\partial}{\partial p_j}\right)
\sum_{\eta} \frac{\chi_{\nu}(C_\eta)}{z_{\eta}} p_{\eta}.
\end{eqnarray*}
See e.g. \cite{Zho1, LLZ}.
From this one easily proves the following results.

\begin{proposition} \label{prop:CutJoin}
We have the following equations:
\begin{eqnarray}
&& \frac{\pa \Phi_h^\bu}{\pa\lambda}
= \frac{1}{2} \sum_{i, j\geq 1}
\left(ijp^{\pm}_{i+j}\frac{\pa^2\Phi_h^\bu}{\pa p^\pm_i\pa p^\pm_j}
+ (i+j)p^\pm_i p^\pm_j\frac{\pa \Phi_h^\bu}{\pa p^\pm_{i+j} }\right), \label{eqn:CutJoin3} \\
&& \frac{\pa \Phi_h^{\circ}}{\pa \lambda}
= \frac{1}{2} \sum_{i, j\geq 1} 
  \left(ijp^\pm_{i+j}\frac{\pa^2\Phi_h^{\circ}}{\pa p_i^\pm\pa p_j^\pm}
+ ijp^\pm_{i+j}\frac{\pa\Phi_h^{\circ}}{\pa p^\pm_i}\frac{\pa \Phi_h^{\circ}}{\pa p^\pm_j}
+ (i+j)p^\pm_i p^\pm_j\frac{\pa \Phi_h^{\circ}}{\pa p^\pm_{i+j}}\right).  \label{eqn:CutJoin4}
\end{eqnarray}
\end{proposition}

The new feature for the double Hurwitz numbers is that there are two choices
to do the cut-and-join, on the $+$  side or on the $-$ side.
One can rewrite (\ref{eqn:CutJoin3}) as sequences of systems of ODEs as follows.
For each partition $\mu^-$ of $d$,
one gets a system of ODEs
for $\{\Phi_{\mu^+, \mu^-}^{\bullet}(\lambda): |\mu^+|=d\}$,
hence they are determined by
$\{\Phi^{\bullet}_{\mu^+, \mu^-}(0): |\mu^+|=d\}$.
One can also reverse the roles of $\mu^+$ and $\mu^-$.
There are matrices $CJ_d$ such that the cut-and-join equations in degree $d$ can be written as
\begin{eqnarray} \label{eqn:CJPhid}
&& \frac{d}{d \lambda} \Phi^{\bullet}_d = CJ_d \cdot \Phi^{\bullet}_d
= \Phi^{\bullet}_d \cdot CJ^t_d.
\end{eqnarray}

\begin{ex}
{\em  When $d=2$,
the cut-and-join equation becomes
\begin{eqnarray*}
\frac{d}{d\lambda}\begin{pmatrix} \Phi^{\bullet}_{(2), (2)} & \Phi^{\bullet}_{(2), (1^2)} \\
\Phi^{\bullet}_{(1^2), (2)} & \Phi^{\bullet}_{(1^2), (1^2)} \end{pmatrix}
& = & \begin{pmatrix} 0 & 1 \\ 1 & 0 \end{pmatrix}
\begin{pmatrix} \Phi^{\bullet}_{(2), (2)} & \Phi^{\bullet}_{(2), (1^2)} \\
\Phi^{\bullet}_{(1^2), (2)} & \Phi^{\bullet}_{(1^2), (1^2)} \end{pmatrix} \\
& = & \begin{pmatrix} \Phi^{\bullet}_{(2), (2)} & \Phi^{\bullet}_{(2), (1^2)} \\
\Phi^{\bullet}_{(1^2), (2)} & \Phi^{\bullet}_{(1^2), (1^2)} \end{pmatrix}
\begin{pmatrix} 0 & 1 \\ 1 & 0 \end{pmatrix}
\end{eqnarray*}
The initial values are:
\begin{eqnarray*}
&& \begin{pmatrix} \Phi^{\bullet}_{(2), (2)} & \Phi^{\bullet}_{(2), (1^2)} \\
\Phi^{\bullet}_{(1^2), (2)} & \Phi^{\bullet}_{(1^2), (1^2)} \end{pmatrix}(0)
= \begin{pmatrix} \frac{1}{2} & 0 \\ 0 & \frac{1}{2} \end{pmatrix}
\end{eqnarray*}
Hence we have the following solution:
\begin{eqnarray*}
&& \begin{pmatrix} \Phi^{\bullet}_{(2), (2)} & \Phi^{\bullet}_{(2), (1^2)} \\
\Phi^{\bullet}_{(1^2), (2)} & \Phi^{\bullet}_{(1^2), (1^2)} \end{pmatrix}(\lambda)
= \begin{pmatrix}
\frac{1}{2} \cosh \lambda & \frac{1}{2} \sinh \lambda \\
\frac{1}{2} \sinh \lambda & \frac{1}{2} \cosh \lambda
 \end{pmatrix}
\end{eqnarray*}
This is compatible with (\ref{eqn:BUR2}).}
\end{ex}

\subsection{Cut-and-join equation for two-partition Hodge integrals}
\label{sec:Proof}
For each $(\mm)\in\cP^2_+$, we will define a generating function
$$
K^\bu_\mm(\lam)
$$
of relative Gromov-Witten invariants. In Section \ref{sec:LocalizeX},
we will derive the following identity by relative virtual localization:
\begin{eqnarray*}
\textup{(\ref{eqn:KG})} && K^\bu_\mm(\lam)
= \sum_{|\nu^{\pm}| = |\mu^{\pm}|}
\Phi^\bullet_{\mu^+, \nu^+}(-\sqrt{-1}\tau\lambda)z_{\nu^+}G^\bu_{\nn}(\lam; \tau)
z_{\nu^-} \Phi^\bu_{\nu^-, \mu^-}(\frac{-\sqrt{-1}}{\tau}\lam).
\end{eqnarray*}
In matrix form, one has for $d^+, d^- \geq 0$,
$$K^\bu(\lam)_{d^+,d^-}
= \Phi^\bu(-\sqrt{-1}\tau\lambda)_{d^+} Z_{d^+} G^\bu(\lam; \tau)_{d^+,d^-}
Z_{d^-} \Phi^{\bu}(\frac{-\sqrt{-1}}{\tau}\lam)_{d^-}.
$$
Hence by (\ref{eqn:PhiInverse}) we have
$$G^\bu(\lambda;\tau)_{d^+,d^-}
= \Phi^\bu(\sqrt{-1}\tau\lambda)_{d^+} Z_{d^+} K^\bu(\lambda)_{d^+,d^-}
Z_{d^-} \Phi^\bu(\frac{\sqrt{-1}}{\tau}\lam)_{d^-}.
$$
Taking derivative in $\tau$ on both sides,
one then gets:
$$ \frac{\pa}{\pa \tau} G^\bu(\lam;\tau)_{d^+,d^-}
= \sqrt{-1}\lam \left(CJ_{d^+} \cdot G^\bu(\lam;\tau)_{d^+,d^-}
- \frac{1}{\tau^2}G^\bu(\lam;\tau)_{d^+,d^-} \cdot CJ_{d^-}^t\right).$$
This completes the proof of the cut-and-join equation for
$G^\bu$ and hence the proof of the formula (\ref{eqn:Conj}) of two-partition
Hodge integrals.

\begin{corollary}\label{zeroframing}
We have
\begin{eqnarray*}
&& K^\bu_\mm(\lam)
= \sum_{\eta^\pm}
\frac{\chi_{\eta^+}(C_{\mu^+})}{z_{\mu^+}}
\cdot
\cW_{\eta^+, \eta^-}(e^{\sqrt{-1}\lam})
\cdot  \frac{\chi_{\eta^-}(C_{\mu^-})}{z_{\mu^-}}.
\end{eqnarray*}
\end{corollary}

\begin{proof}
By (\ref{eqn:KG}), (\ref{eqn:Conj}) and (\ref{eqn:BUR1}) we have
\begin{eqnarray*}
&& K^\bu_\mm(\lam) \\
& = & \sum_{|\nu^{\pm}| = |\mu^{\pm}|}
\Phi^\bullet_{\mu^+, \nu^+}(-\sqrt{-1}\tau\lambda)z_{\nu^+}G^\bu_{\nn}(\lam; \tau)
z_{\nu^-} \Phi^\bu_{\nu^-, \mu^-}(\frac{-\sqrt{-1}}{\tau}\lam) \\
& = & \sum_{\nu^\pm, \rho^\pm, \eta^\pm}e^{-\sqrt{-1}f_{\rho^+}(2)\tau\lambda}
\frac{\chi_{\rho^+}(C_{\mu^+})}{z_{\mu^+}} \frac{\chi_{\rho^+}(C_{\nu^+})}{z_{\nu^+}}
 \cdot z_{\nu^+} \\
&& \cdot
\frac{\chi_{\eta^+}(C_{\nu^+})}{z_{\nu^+}}
e^{\sqrt{-1}\kappa_{\eta^+} \tau\lambda/2}
\cW_{\eta^+, \eta^-}(e^{\sqrt{-1}\lam}) e^{\sqrt{-1}\kappa_{\eta^-} \tau^{-1}\lambda/2}
\frac{\chi_{\eta^-}(C_{\nu^-})}{z_{\nu^-}} \\
&& \cdot z_{\nu^-}
\cdot \frac{\chi_{\rho^-}(C_{\nu^-})}{z_{\nu^-}} \frac{\chi_{\rho^-}(C_{\mu^-})}{z_{\mu^-}}
e^{-\sqrt{-1}f_{\rho^-}(2)\tau^{-1}\lambda} \\
& = & \sum_{\eta^\pm}
\frac{\chi_{\eta^+}(C_{\mu^+})}{z_{\mu^+}}
\cdot
\cW_{\eta^+, \eta^-}(e^{\sqrt{-1}\lam})
\cdot  \frac{\chi_{\eta^-}(C_{\mu^-})}{z_{\mu^-}}.
\end{eqnarray*}
In the last equality we have used (\ref{eqn:Orth1}).
\end{proof}

\section{Relative stable morphisms and relative virtual localization}
\label{sec:LiGV}
In this section, we  will give a brief review of the moduli spaces of
algebraic relative stable morphisms
\cite{Li1, Li2} (see \cite{Li-Rua, Ion-Par1, Ion-Par2} for the symplectic
approach) and virtual localization on such spaces \cite{Gra-Pan, Gra-Vak2}.

\subsection{Relative stable morphisms} \label{relative}
The  definitions given in this section are based on J. Li's works on relative stable
morphisms \cite{Li1, Li2}, with minor modifications.

Let $Y$ be a smooth projective variety. Let $D^1,\ldots, D^k$ be
disjoint smooth divisors in $Y$.
For $\alpha=1,\ldots,k$, define
$$
\Da=\bP(\cO_{\da}\oplus \cN_{\da/Y}) \to \da,
$$
where $\cN_{D/Y}$ denotes the normal sheaf of a subvariety
$D$ in $Y$. The projective line bundle
$\Da \to \da$ has two distinct sections
$$
\da_0     =\bP(\cO_{\da}\oplus 0),\ \ \
\da_\infty=\bP(0\oplus \cN_{\da/Y}).
$$
We have
$$
\cN_{\da_0/\Da}\cong \cN_{\da/Y}^{-1},\ \ \
\cN_{\da_\infty /\Da}\cong \cN_{\da/Y}.
$$

Let
$$
\Da(m)=\Da_1\cup \Da_2\cup\cdots\cup \Da_m,
$$
where $\Da_i\cong \Da$ for $i=1,\ldots,m$. Let
$\da_{i,0}$ and $\da_{i,\infty}$ be the two
distinct sections of $\Da_i$ which correspond
to $\da_0$ and $\da_\infty$, respectively. Then
$\Da(m)$ is obtained by identifying
$\da_{i,\infty}$ with $\da_{i+1,0}$ for $i=1,\cdots, m-1$
under the canonical isomorphisms
$$
\da_{i,\infty}\cong \da \cong \da_{i+1,0}.
$$
Define
$$
\da_{(0)}=\da_{1,0},\ \ \
\da_{(i)}=\Da_i\cap \Da_{i+1},\ \ \
\da_{(m)}= \da_{m,\infty},
$$
where $i=1,\ldots, m-1$. The $\bC^*$ action on $\cO_{\da}$ induces
a $\bC^*$ action on $\Da$ such that $\Da\to \da$ is $\bC^*$ equivariant,
where $\bC^*$ acts on $\da$ trivially. The two distinct sections
$\da_0, \da_\infty$ are fixed under this $\bC^*$ action.
So there is a $(\bC^*)^m$ action on $\Da(m)$ fixing
$\da_{(0)},\ldots, \da_{(m)}$, such that $\Da(m)\to \da$ is
$(\bC^*)^m$ equivariant, where $(\bC^*)^m$ acts on $\da$ trivially.

The variety
$$
\Ym=Y\cup \bigcup_{\alpha=1}^k \Da(m^\alpha)
$$
with normal crossing singularities is obtained by identifying
$\da\subset Y$ with $\da_{(0)}\subset \Da$ under the canonical
isomorphism.  There is a morphism
$$
\pim: \Ym \to Y
$$
which contracts $\Da(m^\alpha)$ to $\da$.
The $(\bC^*)^{m^\alpha}$ action on $\Da(m^{\alpha})$ gives
a $(\bC^*)^{m^1+\cdots+m^k}$ on $\Ym$ such that
$\pim$ is  $(\bC^*)^{m^1 +\cdots+ m^k}$ equivariant
with respect to the trivial action on $Y$.

With the above notation, we are now ready to define
relative stable morphisms for $(Y;D^1,\ldots, D^k)$.

\begin{df}\label{connected}
Let $\beta\in H_2(Y,\bZ)$ be a nonzero homology class such that
$$
d^\alpha=\int_\beta c_1(\cO(\da))\geq 0.
$$
Let $\mu^\alpha$ be a partition of $d^\alpha$.
Define
$$
\MgY
$$
to be the moduli space of morphisms
$$
f:\domain\to \Ym
$$
such that
\begin{enumerate}
\item[(1)] $\domain$
 is a connected prestable curve of arithmetic genus $g$
 with $\sum_{\alpha=1}^k l(\mu^\alpha)$ marked points.
\item[(2)] $(\pim\circ f)_*[C]= \beta\in H_2(Y;\bZ)$.
\item[(3)]
$$
f^{-1}(\da_{(m^\alpha)})=\sum_{i=1}^{l(\mu^\alpha)}\mu^\alpha_i x^\alpha_i
$$
 as Cartier divisors. In particular, if $d^\alpha=0$, then
$f^{-1}(\da_{(m^\alpha)})$ is empty.
\item[(4)] The preimage of $\da_{(l)}$ consists of
      nodes of $C$, where $0\leq l\leq m^\alpha -1$. If
      $f(y)\in \da_{(l)}$ and $C_1$ and $C_2$ are
      two irreducible components of $C$ which intersect at $y$,
      then $f|_{C_1}$ and $f|_{C_2}$ have the same contact order
      to $D^\alpha_{(l)}$ at $y$.
\item[(5)] The automorphism group of $f$ is finite.
\end{enumerate}
Two morphisms described above are isomorphic if they differ by an isomorphism
of the domain and an element in $(\bC^*)^{m^1+\cdots+m^k}$ acting
on the target. In particular, this defines the automorphism
group in the stability condition (5) above.
\end{df}

\begin{rem}
In \cite{Li1, Li2}, the number of divisors $k=1$, but the construction
and proofs in \cite{Li1,Li2} show that
$$
\MgY
$$
is a separated, proper Deligne-Mumford stack with a perfect obstruction theory
of virtual dimension
$$
\int_\beta c_1(TY) + (1-g)(\dim Y-3)+
\sum_{\alpha=1}^k (l(\mu^\alpha)-|\mu^\alpha|),
$$
where $TY$ is the tangent bundle of $Y$.
\end{rem}

\begin{df} \label{disconnected}
We define the moduli space $\MY$ similarly,
with (1) replaced by the following (1)$^\bu$, and
one additional condition (6):
\begin{enumerate}
\item[(1)$^\bu$] $\domain$
 is a possibly disconnected prestable curve with
 $\sum_{\alpha=1}^k l(\mu^\alpha)$ marked points. Let
 $C_1,\ldots, C_n$ be the connected components of $C$, and let
 $g_i$ be the arithmetic genus of $C_i$. Then
$$
\sum_{i=1}^n (2-2g_i) =\chi.
$$
\item[(6)] Let $\beta_i=\tf_*[C_i]$, where $C_i$ is a connected  component
of $C$. Then $\beta_i\neq 0$, and
$$
\int_{\beta_i}c_1(\cO(D^\alpha))\geq 0
$$
for $\alpha=1,\ldots,k$.
\end{enumerate}
\end{df}

The moduli space
$$
\MY
$$
is a finite quotient of a disjoint union of products of the moduli spaces
defined in Definition \ref{connected}. By \cite{Li1, Li2}, it is a
separated, proper Deligne-Mumford stack with a perfect obstruction theory
of virtual dimension
$$
\int_\beta c_1(TY) + \frac{\chi}{2}(\dim Y-3)
+\sum_{\alpha=1}^k (l(\mu^\alpha)-|\mu^\alpha|).
$$

\subsection{Tangent and obstruction spaces}
This section is based on \cite[Section 5.1]{Li2}. We first introduce some notation.
If $m^\alpha>0$, define line bundles $L^\alpha_l$ on $\da_{(l)}\subset Y[m^1,\ldots,m^k]$ by 
$$
L^\alpha_l=
\left\{\begin{array}{ll}
N_{\da_{(0)}/Y}\otimes N_{\da_{(0)}/\Da_1} &l=0\\
N_{\da_{(l)}/\Da_l}\otimes N_{\da_{(l)}/\Da_{l+1}} & 1\leq l\leq m^\alpha-1
\end{array}\right.
$$
Note that $L^\alpha_l$ is a trivial line bundle on $\da_{(l)}$. 

The tangent space $T^1$ and the obstruction space $T^2$ of
$$
\MgY
$$
at the moduli point
$$
\left[f: \domain \to \Ym \ \right]
$$
are given by the following two exact sequences:

\begin{equation}\label{eqn:exactI}
\begin{array}{ccc}
0 &\to \Ext^0(\Omega_C(R), \cO_C) \to H^0(\mathbf{D}^\bu)\to  T^1 &    \\
  &\to \Ext^1(\Omega_C(R), \cO_C) \to H^1(\mathbf{D}^\bu)\to  T^2 &\to  0
\end{array}
\end{equation}
\begin{equation}\label{eqn:exactII}
\end{equation}
\begin{eqnarray*}
0 & \to&  H^0\left(C,\TY\right)\to H^0(\mathbf{D}^\bu)
     \to \bigoplus_{\alpha=1}^k \bigoplus_{l=0}^{m^\alpha -1} \zero \\
  & \to& H^1\left(C,\TY\right) \to H^1(\mathbf{D}^\bu)
      \to \bigoplus_{\alpha=1}^k \bigoplus_{l=0}^{m^\alpha -1} \one \to 0
\end{eqnarray*}
where
$$
R=\sum_{\alpha=1}^k\sum_{i=1}^{l(\mu^\alpha)}x_i^\alpha,
$$
\begin{eqnarray}
&&\zero \cong  \bigoplus_{q\in f^{-1}(\da_{(l)}) }
T_q (f^{-1}(\Da_l))\otimes T^*_q (f^{-1}(\Da_l))
\cong \bC^{\oplus n^\alpha_l},\\
&&\one\cong \left. H^0(\da_{(l)},L^\alpha_l)^{\oplus n^\alpha_l}\right/
H^0(\da_{(l)}, L^\alpha_l),\label{eqn:one}
\end{eqnarray}
and $n^\alpha_l$ is the number of nodes over $D^\alpha_l$. In 
(\ref{eqn:one}), 
$$
H^0(\da_{(l)},L^\alpha_l)\to H^0(\da_{(l)},L^\alpha_l)^{\oplus n^\alpha_l}
$$
is the diagonal embedding. 

We refer the reader to  \cite{Li2} for the definitions
of $H^i(\mathrm{D}^\bu)$  and the maps between terms in (\ref{eqn:exactI}), 
(\ref{eqn:exactII}). Here we only explain the part
relevant to virtual localization calculations. The vector space 
$$
B_1= \Ext^0(\Omega_C(R), \cO_C)
$$
is the space of the infinitesimal automorphisms of the domain curve
$(C,R)$, and
$$
B_4= \Ext^1(\Omega_C(R), \cO_C)
$$
is the space of the infinitesimal deformations of $(C,R)$.
Let $\hat{C}$ be the normalization of $C$, $\hat{R}\subset \hat{C}$
be the pull back of $R$, and $R'\subset \hat{C}$ be the divisor
corresponding to nodes in $C$.
From the local to global spectral sequence, we have an exact sequence 
$$
0\to B_{4,0} \to B_4 \to B_{4,1} \to 0,
$$
where
$$
B_{4,0}=H^1(C,\cE xt^0_{\cO_C}(\Omega_C(R),\cO_C))=H^1(C,\Omega_C(R)^\vee)
$$
is the space of infinitesimal deformations of the smooth pointed curve
$(\hat{C},\hat{R} + R')$, and
$$
B_{4,1}=H^0(C,\cE xt^1_{\cO_C}(\Omega_C(R), \cO_C))
\cong \bigoplus_{q\in \mathrm{Sing}(C)}T_{q'} \hat{C}\otimes T_{q''}\hat{C}
$$ 
corresponds to smoothing of nodes of the domain curve. Here
$\mathrm{Sing}(C)$ is the set of nodes of $C$, and 
$q',q''\in \hat{C}$ are the two preimages of $q$ under the normalization map
$\hat{C}\to C$. The tangent line of smoothing of the node $q$ is canonically
identified with  $T_{q'}\hat{C}\otimes T_{q''}\hat{C}$.

The complex vector space
$$
B_2= H^0\left(C,\TY\right)
$$
is the space of infinitesimal deformations of
the map $f$ with fixed domain and target, and
$$
B_5 = H^1\left(C,\TY\right)
$$
is the obstruction space to deforming $f$ with fixed domain and target.

Finally, let
$$
B_3= \bigoplus_{\alpha=1}^k \bigoplus_{l=0}^{m^\alpha -1} \zero,\ \ \
B_6= \bigoplus_{\alpha=1}^k \bigoplus_{l=0}^{m^\alpha -1} \one.
$$
The complex vector space $\one$ correponds to obstruction to smoothing the nodes
in $f^{-1}(\da_{(l)})$. More explicitly, let
$$
f^{-1}(\da_{(l)})=\{ q_1,\ldots, q_n \},
$$
and let $\nu_i$ be the contact order of $f$ to $\da_{(l)}$
at $q_i$ (of either of the two branches of $f$ near $q_i$). Then
$B_4\to H^1(\mathbf{D}^\bu)$ in (\ref{eqn:exactI}) induces a map
\begin{eqnarray*}
\bigoplus_{i=1}^n T_{q_{i,1}}\hat{C}\otimes T_{q_{i,2}}\hat{C}
&\to& \one\cong \left. H^0(\da_{(l)},L_l^\alpha)^{\oplus n}\right/ 
H^0(\da_{(l)},L_l^\alpha)\\
(s_1,\ldots, s_n)& \mapsto&
\left[ (s_1^{\nu_1},\ldots,s_n^{\nu_n}) \right]
\end{eqnarray*}
where we use isomorphisms 
$$
H^0(\da_{(l)},L_l^\alpha)\cong (L_l^\alpha)_{f(q_i)}\cong
 \left(T_{q'_i}\hat{C}\otimes T_{q''_i}\hat{C}\right)^{\otimes \nu_i}.
$$
The first isomorphism follows from the triviality of the line
bundle $L_l^\alpha\to \da_{(l)}$.
We see that the obstruction vanishes iff the smoothing of the nodes
$q_1,\ldots, q_n$ is compatible with the smoothing the 
target along the divisor $\da_{(l)}$, which is parametrized by the 
complex line $H^0(\da_{(l)},L_l^\alpha)$.

\subsection{Relative virtual localization}
In this section, we assume that a torus $T=(\bC^*)^r$ acts on $Y$, and $D^1,\ldots,D^k$
are $T$-invariant divisors. 

Under our assumption, $\cN_{\da/Y}\to \da$ is $T$-equivariant, and the $T$-action
extends to $\Da$. So $T$ acts on $\Ym$, and acts on
$\MgY$ by moving the image. 

The $T$ fixed points set $\MgY^T$ is a disjoint union of
$$
\{ \cF_\Gamma \mid \Gamma\in \GgY \},
$$
where each $\Gamma\in \GgY$ corresponds
to a connected component, or a union of connected components, $\cF_\Gamma$
of $\MgY^T$. Let
$$
[f:\domain\in \cF_\Gamma\subset\MgY,
$$
for some $\Gamma\in G_{g,0}(\bP^1,\mu)$.
The $T$-action on $\MgY$ induces
$T$-actions on the exact sequences (\ref{eqn:exactI}), (\ref{eqn:exactII}) which
define $T^1$ and $T^2$.
Let $T^{i,f}$ and $T^{i,m}$ denote the fixing part
and the moving part of $T^i$ under the $T$-action, respectively,
where $i=1,2$. Then
$$
T^{1,f}- T^{2,f}
$$
defines a perfect obstruction theory
on $\cF_\Gamma$, and
$$
T^{1,m} -T^{2,m}
$$
defines the virtual normal bundle $N^\vir_{\cF_\Gamma}$ of
$\cF_\Gamma$ in $\MgY$. More explicitly, let
$B_i^m$ denote the moving part of $B_i$ under $T$-action, where
$i=1,\ldots,6$. Then $B_3^m=0$. Note that there are subtleties due to the 
$(\bC^*)^{m^1+\cdots + m^k}$ action on the target
$Y[m^1,\ldots,m^k]$. We have 
$$
\frac{1}{e_T(N^\vir_{\cF_\Gamma})}=\frac{e_T(T^{2,m})}{e_T(T^{1,m})}
=\frac{e_T(B_1^m) e_T(B_5^m) e_T(B_6^m)}{e_T(B_2^m) e_T(B_4^m)}
$$

In \cite{Gra-Pan}, T. Graber and R. Pandharipande proved a localization
formula for the virtual fundamental class in the general context of 
$\bC^*$-equivariant perfect obstruction theory.
In \cite{Gra-Vak2}, T. Graber and R. Vakil showed that moduli spaces 
of relative stable morphisms satisfy the technical assumptions required in the
general formalism in \cite{Gra-Pan}, and derived relative virtual
localization under the assumption that the divisor is fixed pointwisely under
the $\bC^*$ action \cite[Theorem 3.6]{Gra-Vak2}. In our context,
the localization formula proved in \cite{Gra-Pan} reads:
\begin{equation}\label{rvl-g}
{[\MgY]}_T^\vir =\sum_{\Gamma\in \GgY} (i_{\cF_\Gamma})_*\left(
\frac{[\cF_\Gamma]_T^\vir}{e_T(N_{\cF_\Gamma}^\vir)}\right)
\end{equation}
where
$$
i_{\cF_\Gamma}:\cF_\Gamma\to \MgY
$$
is the inclusion, $e_T(N^\vir_{\cF_\Gamma})$ is the $T$-equivariant
Euler class of the virtual normal bundle
$N^\vir_{\cF_\Gamma}=T^{1,m}-T^{2,m}$ over $\cF_\Gamma$,
$$
[\MgY]_T^\vir \in A_*^T(\MgY;\bQ)
$$
is the $T$-equivariant virtual fundamental class defined by
the $T$-equivariant perfect obstruction theory $T^1 - T^2$ on $\MgY$, and
$$
[\cF_\Gamma]_T^\vir \in A_*^T (\cF_\Gamma;\bQ)
$$
is the $T$-equivariant virtual fundamental class defined by
the perfect obstruction theory $T^{1,f} - T^{2,f}$ on
$\cF_\Gamma$.

Similarly, we have
\begin{equation}\label{rvl}
{[\MY]}^\vir =\sum_{\Gamma\in \GY} (i_{\cF_\Gamma})_*\left(
\frac{[\cF_\Gamma]^\vir}{e_T(N_{\cF_\Gamma}^\vir)}\right)
\end{equation}

\section{Double Hurwitz numbers as relative Gromov-Witten invariants}
\label{sec:LocalizeP}

In this section, we study double Hurwitz numbers by
relative Gromov-Witten theory.

\subsection{Relative morphisms to $\bP^1$} \label{line}
Let $[Z_0,Z_1]$ be the homogeneous coordinates of
$\bP^1$. Let $\bC^*$ act on $\bP^1$ by
$$
t\cdot[Z_0,Z_1]=[t Z_0, Z_1]
$$
for $t\in \bC^*$, $[Z_0, Z_1]\in \bP^1$.
Let
$$
s^+ = [0,1],\ \ \ s^- = [1,0]
$$
be the two fixed points of this $\bC^*$-action.

Let $\mu^+$ and $\mu^-$ be two partitions of $d>0$.
Let $[\bP^1]\in H_2(\bP^1;\bZ)$ be the fundamental class.
Define
\begin{eqnarray*}
\MgP&=&\Mbar_{g,0}(\bP^1,s^+,s^-; d[\bP^1],\mu^+,\mu^-),\\
\MP&=&\Mbar^\bullet_\chi(\bP^1,s^+,s^-; d[\bP^1],\mu^+,\mu^-).\\
\end{eqnarray*}
The virtual dimension of $\MgP$ is
$$
\con,
$$
and the virtual dimension of $\MP$ is
$$
\dis.
$$

We extend the $\bC^*$ action on $\bP^1$ to
$\bP^1[m^+,m^-]$ by trivial action on
$\Delta^\pm[m^\pm]$, which is a chain of
$m^\pm$ copies of $\bP^1$. This induces
$\bC^*$-actions on $\MgP$ and $\MP$.
Define the moduli spaces of unparametrized relative stable maps
to the triple $(\bP^1,s^+,s^-)$ to be
\begin{eqnarray*}
\MgP//\bC^* &=&\left.\left(\MgP\setminus\MgP^{\bC^*}\right)\right/\bC^*, \\
\MP //\bC^* &=&\left.\left(\MP \setminus\MP ^{\bC^*}\right)\right/\bC^*.
\end{eqnarray*}

Then $\MgP//\bC^*$ is a separated, proper Deligne-Mumford stack with a
perfect obstruction theory of virtual dimension
$$
\con-1,
$$
and $\MP//\bC^*$ is a separated, proper Deligne-Mumford stack with a
perfect obstruction theory of virtual dimension
$$
\dis-1.
$$

\subsection{Target $\psi$ classes}
In the notation in Section \ref{relative}, we have
$\Delta^\pm\cong\bP^1$, $\Delta^\pm(m)$ is a chain of $m$ copies of $\bP^1$,
and $D^\pm_{(l)}$ is a point, for $l=0,\ldots,m^\pm$.  Let
$\bL^\pm$ and  be the line bundle on $\MgP//\bC^*$ whose fiber at
$$
\left[ f:(C,x_1,\ldots,x_{l(\mu^+)}, y_1,\ldots,y_{l(\mu^-)})
\to \bP^1[m^+,m^-] \right]\in \MgP
$$
is the cotangent line
$$
T^*_{D^\pm_{(m^\pm)}}\left(\bP^1[m^+,m^-]\right)
$$
of $\bP^1[m^+,m^-]$ at the smooth point $D^\pm_{(m^\pm)}$.
We define $\bL^+$ and $\bL^-$  on $\MP//\bC^*$,
similarly. Define the {\em target $\psi$ classes}
$$
\psi^0 = c_1(\bL^+),\ \ \ \psi^\infty=c_1(\bL^-).
$$

The following integral of $\psi^0$ arises in the localization calculations in \cite{LLZ}:
$$
\int_{[\MgP//\bC^*]^\vir}(\psi^0)^{\con-1}.
$$
In Section \ref{sec:ELSV}, we will relate such integrals of target $\psi$
classes to double Hurwitz numbers (Proposition \ref{pro:xELSV}, \ref{pro:gELSV}).

\subsection{Double Hurwitz numbers}
Let $\mm$ be two partitions of $d>0$.
There are branch morphisms
\begin{eqnarray*}
\Br:\MgP &\to& \Sym^\con\bP^1\cong\bP^\con\\
\Br: \MP &\to& \Sym^\dis\bP^1\cong \bP^\dis
\end{eqnarray*}
The double Hurwitz numbers for connected covers of $\bP^1$ can
be defined by
$$
H^\circ_g(\mm)=\frac{1}{|\Aut(\mu^+)||\Aut(\mu^-)|}\int_{[\MgP]^{\vir} }\Br^*(H^{\con})
$$
where $H\in H^2(\bP^\con;\bZ)$ is the hyperplane
class. The double Hurwitz numbers for possibly disconnected covers
of $\bP^1$ can be defined by
$$
\xm{H} =\frac{1}{\amm}\int_{[\MP]^{\vir} }\Br^*( H^{\dis}).
$$
We have
$$
H^\bu_{2-2g,\mm}=H^\bu_g(\mm).
$$
Recall that  $H^\circ_g(\mm), H^\bu_g(\mm)$ are defined combinatorially in
Section \ref{sec:Hurwitz}.

We define generating functions of double Hurwitz numbers as in Section 3:
\begin{eqnarray*}
\Phi^\circ_\mm(\lam)&=&\sum_{g=0}^\infty\frac{\lam^\con}{(\con)!}H_g^\circ(\mm)\\
\Phi_\mm^\bu(\lam) &=&
\sum_{\chi\in 2\bZ, \frac{\chi}{2}\leq \min\{l(\mu^+),l(\mu^-)\} }
\frac{\lam^\dis}{(\dis)!}\xm{H}\\
\Phi^\circ(\lam;p^+,p^-)&=&\sum_\mm\Phi^\circ_\mm(\lam)p^+_{\mu^+} p^-_{\mu^-}\\
\Phi^\bu(\lam;p^+,p^-)&=&1+ \sum_\mm\Phi_\mm^\bu(\lam)
                    p^+_{\mu^+} p^-_{\mu^-}
\end{eqnarray*}
Then
$$
\Phi^\bu(\lam;p^+,p^-) =\exp(\Phi^\circ(\lam;p^+,p^-)).
$$

Note that
\begin{equation}\label{double-zero}
\Phi_\mm^\bu(0)=\frac{\delta_\mm}{z_{\mu^+}},
\end{equation}
where $z_\nu = \nu_1\cdots\nu_{l(\nu)}|\Aut(\nu)|$,
so
$$
\Phi^\bu(0,p^+,p^-)=1+\sum_{\mu}\frac{p^+_\mu p^-_\mu}{z_\mu}.
$$

\subsection{Gluing formula}\label{sec:glue}

Let $k^+, k^-$ be positive integers such that
$$
k^+ + k^- =\dis.
$$
By gluing formula of algebraic relative Gromov-Witten invariants \cite[Corollary 3.16]{Li2},
we have
\begin{eqnarray*}
&&\int_{[\MP]^\vir}\Br^*(H^\dis)\\
&=&\sum_{-\chi^\pm+l(\mu^\pm)+l(\nu)=k^\pm}\int_{[\Mmn]^\vir} \Br^*(H^{k^+})\\
&&\cdot\frac{a_\nu}{|\Aut(\nu)|}\int_{[\Mnm]^\vir} \Br^*(H^{k^-})
\end{eqnarray*}
where
$$
a_\nu=\nu_1\cdots\nu_{l(\nu)}.
$$

Therefore, we have the following gluing formula for double Hurwitz numbers:
\begin{proposition}[gluing formula]\label{pro:glue}
Let $k^+, k^-$ be positive integers such that
$$
k^+ + k^- = \dis.
$$
Then
\begin{equation}\label{eqn:glue}
\xm{H}=\sum_{-\chi^\pm+l(\mu^\pm)+l(\nu)=k^\pm}
H^\bu_{\chi^+,\mu^+,\nu}z_\nu H^\bu_{\chi^-,\nu,\mu^-}
\end{equation}
\end{proposition}
Recall that $z_\nu=a_\nu|\Aut(\nu)|$.

Let $d=|\mu^+|=|\mu^-|$.
It is straightforward to check that Proposition \ref{pro:glue} implies
the {\em sum formula}
\begin{equation}\label{eqn:sum}
\sum_{|\nu|=d}\Phi^\bu_\mn(\lam_1)z_\nu\Phi^\bu_\nm(\lam_2)
=\Phi^\bu_\mm(\lam_1 + \lam_2)
\end{equation}
which was derived in Section \ref{sec:sum} from the combinatoric definition.

The cut-and-join equations (\ref{eqn:CJPhid})
for double Hurwitz numbers are
special cases $k_+=1$, $k_-=1$ of Proposition \ref{pro:glue}.
More precisely, differentiate (\ref{eqn:sum}) with repect to $\lam_1$,
and then set $\lam_1=0$. We obtain a cut-and-join equation:
\begin{equation}\label{plus-cut-join}
\frac{d}{d\lam}\Phi^\bu_\mm(\lam)
=\sum_{|\nu|=d}H^\bu_{l(\mu^+)+l(\nu)-1, \mn}z_\nu \Phi^\bu_\nm(\lam).
\end{equation}
Differentiate (\ref{eqn:sum}) with repect to $\lam_2$,
and then set $\lam_2=0$. We obtain another cut-and-join equation:
\begin{equation}\label{minus-cut-join}
\frac{d}{d\lam}\Phi^\bu_\mm(\lam)
=\sum_{|\nu|=d}\Phi^\bu_\mn(\lam)z_\nu H^\bu_{l(\nu)+l(\mu^-)-1, \nm}.
\end{equation}

Define the {\em cut-and-join coefficients}
$$
(CJ)_{\mu\nu}=H^\bu_{l(\mu)+l(\nu)-1,\mu,\nu}z_\nu.
$$
They are the entries of the matrix $CJ_d$ in Section  \ref{sec:Hcj}.
The cut-and-join equations can be written as
\begin{equation}
\frac{d}{d\lam}\Phi^\bu_\mm(\lam)
=\sum_{|\nu|=d} (CJ)_{\mu^+\nu}\Phi^\bu_\nm(\lam)
=\sum_{|\nu|=d} \Phi^\bu_\mn(\lam) (CJ)_{\mu^-\nu},
\end{equation}
which is equivalent to (\ref{eqn:CJPhid}) in Section \ref{sec:Hcj}:
$$
\frac{d}{d\lam}\Phi^\bu_d= CJ_d \cdot\Phi^\bu_d =\Phi^\bu_d\cdot CJ^t_d.
$$

\begin{rem}
The cut-and-join equation of Hurwitz numbers $H^\circ_g(\mu),
H^\bu_g(\mu)$ was first proved using combinatorics by Goulden,
Jackson and Vainstein \cite{Gou-Jac-Vai} and later proved  using
gluing formula of symplectic relative Gromov-Witten invariants by
Li-Zhao-Zheng \cite{Li-Zha-Zhe} and Ionel-Parker \cite{Ion-Par2}.
The gluing formula of symplectic relative Gromov-Witten invariants
was proved in \cite{Li-Rua, Ion-Par2}. 
\end{rem}

\subsection{Localization}\label{sec:branch}
In the spirit of \cite[Section 7]{LLZ}, we lift
$$H^\dis\in H^{2(\dis)}(\bP^\dis;\bZ)$$
to
$$
\prod_{k=1}^{\dis}(H-w_ku)\in H^{2(\dis)}_{\bC^*}(\bP^\dis;\bZ),
$$
where $w_k\in \bZ$, and  compute
$$
\xm{H} =\frac{1}{\amm}\int_{[\MP]^{\vir} }\Br^*\left(\prod_{k=1}^\dis(H-w_k u) \right)
$$
by virtual localization.

\subsection{Torus fixed points and admissible triples}
Given a morphism
$$
f:\domain \to \bP^1[m^+,m^-]
$$
which represents a point in $\MP^{\bC^*}$,
let
$$
\tf=\pi[m^+,m^-]\circ f:C\to \bP^1,
$$
and let $C^\pm=\tf^{-1}(s^\pm)$. Then
$$
C=C^+\cup L \cup C^-,
$$
where $L$ is a disjoint union of projective lines.
Let
\begin{eqnarray*}
f^\pm=f|_{C^\pm}: C^\pm &\to& \dpm(m^\pm),\\
f^0= f|_L: L &\to& \bP^1.
\end{eqnarray*}
Then $f^0$ is a morphism of degree
$$
d=|\mu^+|=|\mu^-|
$$
fully ramified over $s^+$ and $s^-$. The degrees of $f^0$ restricted
to connected components of $L$ determine a partition $\nu$ of $d$.

Let $C^+_1,\ldots, C^+_k$ be the connected components of $C^+$, and let
$g_i$ be the arithmetic genus of $C^+_i$. (We define $g_i=0$ if
$C^+_i$ is a point.) Define
$$
\chi^+ =\sum_{i=1}^k(2-2g_i),
$$
and define $\chi^-$ similarly. We have
$$
-\chi^+ + 2 l(\nu) -\chi^- = -\chi.
$$
Note that $\chi^\pm\leq 2\min\{l(\mu^\pm),l(\nu)\}$. So
$$
\xmn\geq 0,
$$
and the equality holds if and only if $m^+=0$. In this case, we have
$\nu=\mu^+$, $\chi^+= 2l(\mu^+)$, and $\chi^-=\chi$. Similarly,
$$
\xnm\geq 0,
$$
and the equality holds if and only if $m^-=0$. In this case, we have
$\nu=\mu^-$, $\chi^-= 2l(\mu^-)$, and $\chi^+=\chi$.
There are three cases:

\paragraph{Case 1: $m^-=0$.} Then $f^-$ is a constant map,
$\chi^+=\chi$, $\nu=\mu^-$, and $f^+$ represents a point in
$$
\MP//\bC^*.
$$
\paragraph{Case 2: $m^+=0$.} Then $f^+$ is a constant map,
$\chi^-=\chi$, $\nu=\mu^+$, and $f^-$ represents a point in
$$
\MP//\bC^*.
$$
\paragraph{Case 3: $m^+, m^- >0$.} Up to an element of
$\Aut(\nu)$, $f^+$ represents
a point in
$$
\Mmn//\bC^*,
$$
and $f^-$ represents an element
of
$$
\Mnm//\bC^*.
$$

\begin{df}
We say a triple $\tri$ is {\em admissible} if
\begin{itemize}
\item $\chi^+,\chi^-\in 2\bZ$.
\item $\nu$ is a partition of $d$.
\item $\chi^\pm\leq 2\min\{l(\mu^\pm), l(\nu)\}$.
\item $-\chi^+ + 2 l(\nu) -\chi^- =  -\chi$.
\end{itemize}
Let $\GP$ denote the set of all admissible triples.
\end{df}

We define
\begin{eqnarray*}
\xmm{\Mbar} &=& \MP//\bC^*,\\
\mmx{\Mbar} &=& \MP//\bC^*,\\
\end{eqnarray*}
and define
$$
\xnx{\Mbar}=
\left(\Mmn//\bC^*\right)\times\left(\Mnm//\bC^*\right).
$$
if $\tri\in\GP$, and
$$
\xmn >  0,\ \ \ \xnm>0.
$$

For every
$\tri\in \GP$, there is a morphism
$$
\xnx{i}:\xnx{\Mbar}\to \MP,
$$
whose image $\xnx{\cF}$ is a union of connected components of $\MP^{\bC^*}$.
The morphism $\xnx{i}$ induces an isomorphism
$$
\xnx{\Mbar}/\xnx{A} \cong \xnx{\cF},
$$
where
$$
\xmm{A}=\prod_{i=1}^{l(\mu^-)}\bZ_{\mu^-_i} ,\ \ \
\mmx{A}=\prod_{i=1}^{l(\mu^+)}\bZ_{\mu^+_i},
$$
and for $-\chi^\pm+l(\mu^\pm)+l(\nu)>0$, we have
$$
1\to \prod_{i=1}^{l(\nu)}\bZ_{\nu_i}\to \xnx{A}\to\Aut(\nu) \to 1.
$$

Recall that $a_\nu=\nu_1\cdots\nu_{l(\nu)}$, and $z_\nu= a_\nu|\Aut(\nu)|$. We have
$$
|\xmm{A}|=a_{\mu^-},\ \ \
|\mmx{A}|=a_{\mu^+},
$$
and
$$
|\xnx{A}|=z_\nu
$$
if $-\chi^\pm+l(\mu^\pm)+l(\nu)>0$.

The fixed points set $\MP^{\bC^*}$ is a disjoint union of
$$
\{ \xnx{\cF} \mid \tri\in\GP \}
$$

\subsection{Contribution from each admissible triple}
Let $\tri\in\GP$. We have
\begin{eqnarray*}
\Br(\xnx{\cF})&=&(\xmn )s^+ + (\xnm)s^- \\
              & &\in \Sym^\dis \bP^1 =\bP^\dis,
\end{eqnarray*}
so
\begin{eqnarray*}
& &\xnx{i}^*\Br^*\left(\prod_{l=1}^\dis(H-w_l)\right)\\
&=&\left(\prod_{l=1}^\dis(-\chi^+ +l(\mu^+)+l(\nu)-w_l)\right)u^\dis.
\end{eqnarray*}

Let $\xnx{N}^\vir$ on $\xnx{\Mbar}$ be the pull-back of the
virtual normal bundle of $\xnx{\cF}$ in $\MP$. Calculations
similar to those in \cite[Appendix A]{LLZ} show that
\begin{eqnarray*}
\frac{1}{e_{\bC^*}(\xmm{N}^{\vir})}
&=&\frac{a_{\mu^-}}{u-\psi^\infty},\\
\frac{1}{e_{\bC^*}(\mmx{N}^{\vir})}
&=&\frac{a_{\mu^+}}{-u-\psi^0},
\end{eqnarray*}
and for $-\chi^\pm+l(\mu^\pm)+l(\nu)>0$, we have
$$
\frac{1}{e_{\bC^*}(\xnx{N}^{\vir})}
=\frac{a_\nu}{u-\psi_+^\infty}\frac{a_\nu}{-u-\psi_-^0},
$$
where $\psi_+^\infty$, $\psi_-^0$ are the target $\psi$ classes on
$$
\Mmn,\ \ \ \Mnm,
$$
respectively.

Let $w=(w_1,\ldots,w_l)$. Then
\begin{eqnarray*}
&&\xmm{I}(w)\\
&=&\frac{1}{a_{\mu^-}}\int_{[\xmm{\Mbar}]^{\vir}}
\frac{\xmm{i}^*\Br^*\left(\prod_{l=1}^\dis(H-w_l u)\right)}
{e_{\bC^*}(\xmm{N}^{\vir})}\\
&=&\left(\prod_{l=1}^\dis(\dis-w_l)\right)\int_{[\xmm{\Mbar}]^{\vir}}
\frac{u^\dis}{u-\psi^\infty}\\
&=&\left(\prod_{l=1}^\dis(\dis-w_l)\right)\int_{[\MP//\bC^*]^{\vir}}
(\psi^\infty)^{\dis-1}
\end{eqnarray*}
\begin{eqnarray*}
&&\mmx{I}(w)\\
&=&\frac{1}{a_{\mu^+}}\int_{[\mmx{\Mbar}]^{\vir}}
\frac{\mmx{i}^*\Br^*\left(\prod_{l=1}^\dis(H-w_l u)\right)}
{e_{\bC^*}(\mmx{N}^{\vir})}\\
&=&\left(\prod_{l=1}^\dis(-w_l)\right)\int_{[\mmx{\Mbar}]^{\vir}}
\frac{u^\dis}{-u-\psi^0}\\
&=&\left(\prod_{l=1}^\dis w_l\right)\int_{[\MP//\bC^*]^{\vir}}
(\psi^0)^{\dis-1}
\end{eqnarray*}
\begin{eqnarray*}
&&\xnx{I}(w)\\
&=& \frac{1}{z_\nu}\int_{[\xnx{\Mbar}]^{\vir}}
    \frac{\xnx{i}^*\Br^*\left(\prod_{l=1}^\dis(H-w_l u)\right)}
    {e_{\bC^*}(\xnx{N}^{\vir})}\\
&=& \frac{a_\nu}{|\Aut(\nu)|}\left(\prod_{l=1}^\dis(\xmn-w_l)\right)\int_{[\xnx{\Mbar}]^{\vir}}
    \frac{u^\dis}{(u-\psi_+^\infty)(-u-\psi_-^0)}\\
&=& \frac{a_\nu}{|\Aut(\nu)|}\left(\prod_{l=1}^\dis(\xmn-w_l)\right)(-1)^{\xnm}\\
& & \cdot\int_{[\Mmn//\bC^*]^{\vir}}(\psi^\infty)^{\xmn-1}
         \int_{[\Mnm//\bC^*]^{\vir}}(\psi^0)^{\xnm-1}
 \end{eqnarray*}

\subsection{Sum over admissible triples} \label{sec:ELSV}
We have
\begin{eqnarray*}
\xm{H}&=& \frac{1}{\amm}
      \int_{[\MP]^{\vir}}\Br^* \left(\prod_{l=1}^\dis(H-w_l u)\right)\\
      &=&\frac{1}{\amm} \sum_{\tri\in\GP}\xnx{I}(w).
\end{eqnarray*}
Let $w=(0,1,\ldots, \dis-1)$, we have
\begin{eqnarray*}
\xm{H}&=&\frac{1}{\amm}\xmm{I}(w)\\
&=&\frac{(\dis)!}{\amm}\int_{[\MP//\bC^*]^{\vir} }(\psi^\infty)^{\dis-1}.
\end{eqnarray*}
Let $w=(1,2,\ldots,\dis)$, we have
\begin{eqnarray*}
\xm{H}&=&\frac{1}{\amm}\mmx{I}(w)\\
&=&\frac{(\dis)!}{\amm}\int_{[\MP//\bC^*]^{\vir} }(\psi^0)^{\dis-1}.
\end{eqnarray*}
So we have
\begin{proposition}\label{pro:xELSV}
\begin{eqnarray*}\label{xELSV}
&&\frac{\xm{H}}{(\dis)!}\\
&=&\frac{1}{\amm}\int_{[\MP//\bC^*]^{\vir} }(\psi^0)^{\dis-1}\\
&=&\frac{1}{\amm}\int_{[\MP//\bC^*]^{\vir} }(\psi^\infty)^{\dis-1}
\end{eqnarray*}
\end{proposition}

If we replace $\MP$ by $\MgP$ in Section \ref{sec:branch}, we get
\begin{proposition}\label{pro:gELSV}
\begin{eqnarray*}
&&\frac{H^\circ_g(\mm)}{(\con)!}\\
&=&\frac{1}{\amm}\int_{[\MgP//\bC^*]^{\vir} }(\psi^0)^{\dis-1}\\
&=&\frac{1}{\amm}\int_{[\MgP//\bC^*]^{\vir} }(\psi^\infty)^{\dis-1}
\end{eqnarray*}
\end{proposition}

\bigskip

Let $w=(0,1,\ldots, k-1, k+1,\ldots,\dis)$,
where
$$1\leq k\leq \dis-1,$$
we have
\begin{eqnarray*}
&&\xm{H}\\
&=&\frac{1}{\amm}\sum_{\scriptsize \begin{array}{c}\tri\in \GP\\ \xmn=k\end{array}}\xnx{I}(w)\\
&=&\sum_{\scriptsize \begin{array}{c}\tri\in\GP\\ \xmn=k\end{array}}
\frac{(\xmn)!}{|\Aut(\mu^+)|} \int_{[\Mmn]^{\vir} }(\psi^\infty)^{\xmn-1}\\
&&\cdot \frac{a_\nu}{|\Aut(\nu)|}\cdot\frac{(\xnm)!}{|\Aut(\mu^-)|}
\int_{[\Mnm]^{\vir} }(\psi^0)^{\xnm-1}\\
&=&\sum_{\scriptsize \begin{array}{c}\tri\in \GP\\ \xmn=k\end{array}} H^\bu_{\chi^+,\mn}
    z_\nu  H^\bu_{\chi^-,\nm}.
\end{eqnarray*}
This gives an alternative derivation of the gluing formula (\ref{eqn:glue}),
and in particular, the cut-and-join equations (\ref{plus-cut-join}), (\ref{minus-cut-join}).

\section{Moduli spaces and obstruction bundles}
\label{sec:Objects}

In this section, we introduce the geometric objects involved
in the proof of (\ref{eqn:KG}), and fix notation.
\subsection{The target $X$}\label{target}
Let $X$ be the toric surface defined by the fan in Figure 1.
\begin{figure}

\begin{center}
\psfrag{r1}{$\rho_1$}
\psfrag{r2}{$\rho_2$}
\psfrag{r3}{$\rho_3$}
\psfrag{r4}{$\rho_4$}
\psfrag{r5}{$\rho_5$}
\includegraphics{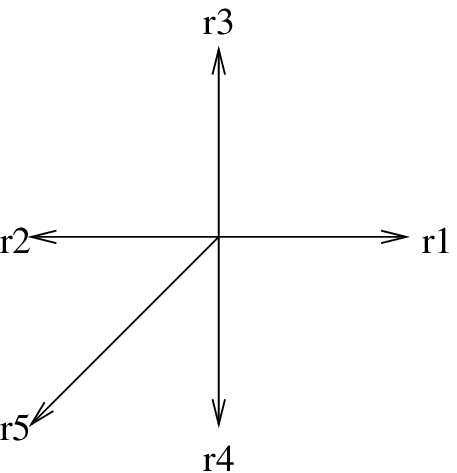}
\end{center}
\caption{The fan of $X$}
\end{figure}
Let $\Phi_i$ be the homogeneous coordinate associated to
the ray $\rho_i$, $i=1, \ldots, 5$, and set
\begin{eqnarray*}
Z_{ij}&=&\{\five{\Phi}\in\bC^5\mid
          \Phi_i=\Phi_j=0\},\\
Z&=& Z_{12}\cup Z_{35}\cup Z_{24}\cup Z_{15}\cup Z_{34}.
\end{eqnarray*}
Then
$$
X=(\bC^5\setminus Z)/(\bC^*)^3,
$$
where $(\bC^*)^3$ acts on $\bC^5$ by
$$
\three{u} \cdot \five{\Phi}=
(u_1\Phi_1,u_1 u_3\Phi_2, u_2\Phi_3, u_2 u_3 \Phi_4, u_3^{-1}\Phi_5),
$$
for $\three{u}\in (\bC^*)^3$, $\five{\Phi}\in\bC^5$.

$T=(\bC^*)^2$ acts on $X$ by
$$
 \two{t} \cdot \Five{\Phi}
=[t_1\Phi_1,\Phi_2,t_2 \Phi_3,\Phi_4,\Phi_5]
$$
for $\two{t}\in T$, $\Five{\Phi}\in X$.

Let
$$
D_i=\{\Five{\Phi}\in X \mid \Phi_i=0 \}
\subset X
$$
be the $T$-invariant divisor associated to the ray $\rho_i$.
Let $\beta_i\in H_2(X;\bZ)$ be the homology class represented
by $D_i$.
We have
\begin{eqnarray*}
H_2(X;\bZ)&=&\left(\bigoplus_{i=1}^5 \bZ \beta_i\right)/
\left(\bZ(\beta_1-\beta_2-\beta_5)\oplus \bZ(\beta_3-\beta_4-\beta_5)
\right)\\
&=& \bZ \beta_1 \oplus \bZ\beta_3 \oplus \bZ\beta_5
\end{eqnarray*}
Let $\beta_i^*\in H^2(X;\bZ)$ be the Poincare dual of $\beta_i$,
$i=1,\ldots,5$.
The intersection form on
$$
H^2(X;\bZ)=\bZ\beta_1^* \oplus \bZ\beta_3^* \oplus \bZ\beta_5^*
$$
is given by
$$
\begin{array}{cccc}
& \beta_1^* & \beta_3^* & \beta_5^* \\
\beta_1^* & 0   & 1   &  0  \\
\beta_3^* & 1   & 0   &  0  \\
\beta_5^* & 0   & 0   & -1
\end{array}
$$
So
$$
 \beta_2^* \cdot \beta_2^*
=\beta_4^* \cdot \beta_4^* = -1,\ \ \
 \beta_2^* \cdot \beta_4^* = 0.
$$

Note that $X$ is a toric blowup of $\bP^1\times\bP^1$ at a point, and
$D_5$ is the exceptional divisor. More explicitly, we have
\begin{eqnarray*}
h: X&\to&\bP^1\times\bP^1\\
\Five{\Phi} &\mapsto&
([\Phi_1,\Phi_2\Phi_5], [\Phi_3,\Phi_4\Phi_5])
\end{eqnarray*}
which is an isomorphism outside $D_5$, and
$h(D_5)=\{([1,0],[1,0])\}$.

The $T$-invariant divisor
$$
K_X=-D_1-D_2-D_3-D_4-D_5
$$
is a canonical divisor of $X$, so
$$
c_1(T_X)=2\beta_1^* + 2\beta_3^* -\beta_5^*.
$$

For $(\mu^+,\mu^-)\in\cP^2_+$, define
$$
\MgX =
\Mbar_{g,0}(X;D_2,D_4\mid |\mu^+|\beta_3 + |\mu^-|\beta_1;\mu^+,\mu^-),
$$
and let $\MX$  be the subset of
\begin{eqnarray*}
\Mbar_\chi^\bullet(X;D_2,D_4\mid  |\mu^+|\beta_3 + |\mu^-|\beta_1;\mu^+,\mu^-)
\end{eqnarray*}
which consists of morphisms
$$
f:C\to X[m^+,m^-]
$$
such that for each connected component $C_i$ of $C$,
$\tf_*[C_i]\in H_2(X;\bZ)$ is an element of
$$
\{ a\beta_3 + b\beta_1 \mid a,b\in \bZ_{\geq 0}, (a,b)\neq (0,0)\}.
$$

The virtual dimension of $\MgX$ is
$$
\gm{r}=g-1 +|\mu^+|+l(\mu^+)+|\mu^-|+l(\mu^-),
$$
and the virtual dimension of $\MX$ is
$$
\xm{r}=-\frac{\chi}{2} +|\mu^+|+l(\mu^+)+|\mu^-|+l(\mu^-).
$$

The moduli space $\MgX$ plays the role of
$\Mbar_{g,0}(\bP^1,\mu)$ in the proof of Mari\~{n}o-Vafa
formula \cite{LLZ}.

We have
$$
D_2\cong\bP^1\cong D_4,\ \ \
\cN_{D_2/X}\cong\cO_{\bP^1}(-1)\cong \cN_{D_4/X},
$$
so
$$
\Delta(D_2)\cong \bF_1\cong \Delta(D_4)
$$
in the notation of Section \ref{relative},
where $\bF_1$ is the Hirzebruch surface
$$
\bP(\cO_{\bP^1}\oplus \cO_{\bP^1}(-1)) \to \bP^1.
$$

\subsection{The obstruction bundles}\label{obstruction}

Let
$$
\pi:\gmu{U}\to\MgX
$$
be the universal domain curve, and let
$$
  P:\gmu{T}\to\MgX
$$
be the universal target.
There is an evaluation map
$$
F:\gmu{U}\to \gmu{T}
$$
and a contraction map
$$
\tpi:\gmu{T}\to X.
$$
Let $\gmu{D}\subset \gmu{U}$ be the
divisor corresponding to the $l(\mu^+)+l(\mu^-)$ marked points.
Define
$$
\gm{V}=R^1\pi_*\left( \tF^*\cO_X(-D_1-D_3)\otimes
\cO_{\gmu{U}}(-\gmu{D})\right)
$$
where $\tF=\tpi\circ F:\gmu{U}\to X$.
The fibers of $\gm{V}$ at
$$
\left[ f: (C,x_1,\ldots,x_{l(\mu^+)}, y_1,\ldots,y_{l(\mu^-)})
      \to X[m^+,m^-] \right]\in \MgX
$$
is
$$
H^1(C, \tf^*\cO_X(-D_1-D_3)\otimes \cO_C(-R))
$$
where $\tf=\pi[m^+,m^-]\circ f$, and
$$
R=x_1+\ldots+x_{l(\mu^+)}+y_1+\cdots+ y_{l(\mu^-)}.
$$
Note that
$$
H^0(C, \tf^*\cO_X(-D_1-D_3)\otimes \cO_C(-R))=0,
$$
and
$$
\deg\tf^*\cO_X(-D_1-D_3)\otimes \cO_C(-R)
=-|\mu^+|-|\mu^-|-l(\mu^+)-l(\mu^-),
$$
so $\gm{V}\to \MgX$ is a vector bundle of rank
$$
\gm{r}=g-1+|\mu^+|+l(\mu^+)+|\mu^-|+l(\mu^-).
$$

The vector bundle $\gm{V}\to \MgX$ plays the role
of the obstruction bundle $V\to\Mbar_{g,0}(\bP^1,\mu)$
in the proof of Mari\~{n}o-Vafa formula
\cite[Section 4.4]{LLZ}.

Similarly, we define a vector bundle $\xm{V}$
of rank
$$
\xm{r}=-\frac{\chi}{2}+|\mu^+|+l(\mu^+)+|\mu^-|+l(\mu^-)
$$
on $\MX$.

\subsection{Torus action}\label{action}
Recall that $T=(\bC^*)^2$ acts on $X$ by
$$
\two{t} \cdot \Five{\Phi}
=[t_1\Phi_1,\Phi_2,t_2\Phi_3,\Phi_4,\Phi_5]
$$
for $\two{t}\in T$,  $\Five{\Phi}\in X$.

Let $T_\bR=U(1)^2$ be the maximal compact subgroup of $T$
The $T_\bR$-action on $X$ determines a moment map
$$
\mu_{T_\bR}:X\to \ft^*_\bR,
$$
where $\ft^*_\bR\cong \bR^2$
is the dual of the Lie algebra $\ft_\bR$ of $T_\bR$.

We now lift the $T$-action on $X$ to the line bundle
$\cO_X(-D_1-D_3)$ as follows. We only need to specify
the representation of $T$ on the fiber of one fixed point
of the $T$ action. The fixed points of the $T$ action
on $X$ are
\begin{eqnarray*}
z_0 = D_1\cap D_3 &=& [0,1,0,1,1]\\
z_+ = D_3\cap D_2 &=& [1,0,0,1,1]\\
z_- = D_1\cap D_4 &=& [0,1,1,0,1]\\
\tilde{z}_+ = D_2\cap D_5 &=& [1,0,1,1,0]\\
\tilde{z}_- = D_4\cap D_5 &=& [1,1,1,0,0]
\end{eqnarray*}

Figure 2 shows the image of $D_1,\ldots,D_5$ and
the above five fixed points under the moment map
$\mu_{T_\bR}:X\to\ft^*_\bR$.

\begin{figure}
\begin{center}
\psfrag{D1}{$D_1$}
\psfrag{D2}{$D_2$}
\psfrag{D3}{$D_3$}
\psfrag{D4}{$D_4$}
\psfrag{D5}{$D_5$}
\psfrag{z0}{$z_0$}
\psfrag{z+}{$z_+$}
\psfrag{z-}{$z_-$}
\psfrag{tz+}{$\tz_+$}
\psfrag{tz-}{$\tz_-$}
\psfrag{a}{$\alpha$}
\psfrag{b}{$\beta$}
\psfrag{-a}{$-\alpha$}
\psfrag{-b}{$-\beta$}
\psfrag{b-a}{$\beta-\alpha$}
\psfrag{a-b}{$\alpha-\beta$}
\includegraphics{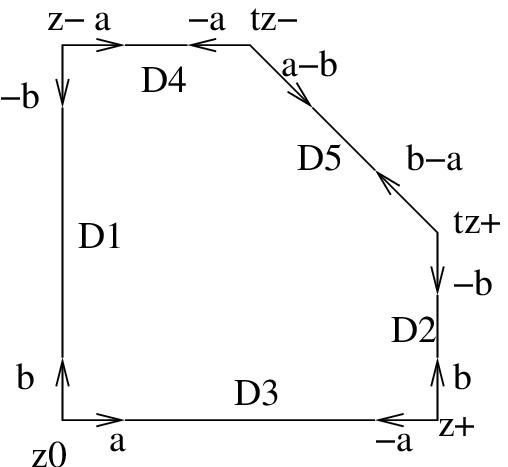}
\end{center}
\caption{The image of $\mu_{T_\bR}:X\to\ft^*_\bR$}
\end{figure}

Let $(w_1,w_2)$ denote the one dimensional representation
given by
$$
\two{t}\cdot z= t_1^{w_1}t_2^{w_2}z
$$
for $\two{t}\in T$, $z\in \bC$. The character ring
of $T$ is given by
$$
\bZ T \cong \bZ[\alpha,\beta],
$$
where $\alpha, \beta$ are the characters of the representations
$(1,0)$, $(0,1)$, respectively. The representations of
$T$ on the fibers of $T_X$ and $\cO_X(-D_1-D_3)$ at fixed
points are given by:
$$
\begin{array}{ccc}
      &        T_X            & \cO_X(-D_1-D_3) \\
  z_0 &   \alpha,\beta        &  -\alpha-\beta  \\
  z_+ & -\alpha,\beta         &   -\beta        \\
  z_- &   \alpha, -\beta      &   -\alpha       \\
\tz_+ & \beta-\alpha, -\beta  &     0           \\
\tz_- & \alpha-\beta, -\alpha &     0
\end{array}
$$

Note that
$$
\ft_\bR^*\cong \bR\alpha\oplus\bR\beta,
$$
and the representations
of $T$ on the fibers of $T_X$ at the fixed points can
be read off from the image of the moment map as
in Figure 3.

The action of $T$ on $\Delta(D_2)$ and $\Delta(D_4)$
can be read off from Figure 3. This extends the action
of $T$ on $X$ to $X[m^+,m^-]$. So $T$ acts on 
$\MgX$ and $\MX$ by moving the image of the morphism.
\begin{figure}
\begin{center}
\psfrag{D2}{$D_{2,0}$}
\psfrag{d2}{$D_{2,\infty}$}
\psfrag{D4}{$D_{4,0}$}
\psfrag{d4}{$D_{4,\infty}$}
\psfrag{a}{$\alpha$}
\psfrag{b}{$\beta$}
\psfrag{-a}{$-\alpha$}
\psfrag{-b}{$-\beta$}
\psfrag{b-a}{$\beta-\alpha$}
\psfrag{a-b}{$\alpha-\beta$}
\includegraphics{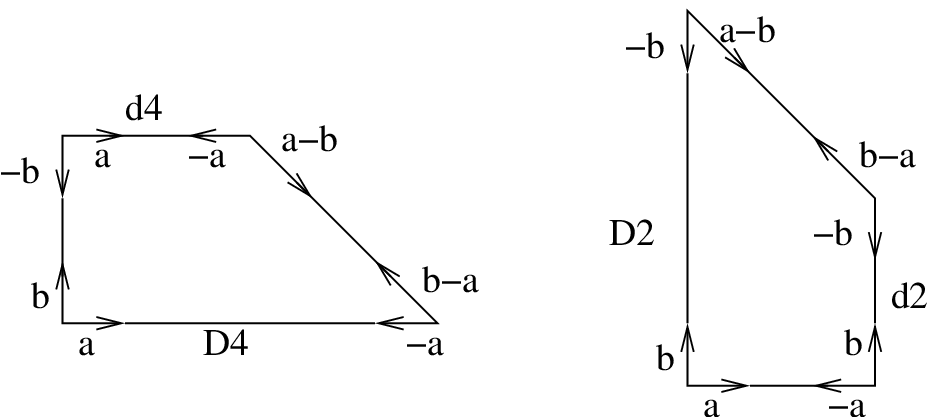}
\end{center}
\caption{The images of $\mu^-_{T_\bR}:\Delta(D_4)\to\ft^*_\bR$
and $\mu^+_{T_\bR}:\Delta(D_2) \to \ft^*_\bR$ }
\end{figure}

The $T$ action on $\cO_X(-D_1-D_3)$ induces $T$ actions
on $\gm{V}$ and on $\xm{V}$. 

\section{Proof of (\ref{eqn:KG})}
\label{sec:LocalizeX}

Let $X$ be defined as in Section \ref{target}. Recall that
$T=(\bC^*)^2$ acts on $X$. Let $D_1,\ldots, D_5$ be $T$-invariant
divisors in $X$ defined in Section \ref{target}.

Let
$$
\xm{V}\to\MX
$$
be defined as in Section \ref{obstruction}, with the torus action defined
in Section \ref{action}. Define
\begin{eqnarray*}
\xm{K}&=&\frac{1}{\amm}\int_{[\MX]^\vir}e(\xm{V}),\\
K^\bu_\mm(\lam)&=&\sum_{\chi\in 2\bZ,\chi\leq 2(l(\mu^+)+l(\mu^-)}\lam^\dis
\frac{(-1)^{|\mu^+|+|\mu^-|} }{\sqrt{-1}^{l(\mu^+)+l(\mu^-)}} \xm{K}.
\end{eqnarray*}

In this section, we will compute
$$
\xm{K}=\frac{1}{\amm}\int_{[\MX]^\vir}e_T(\xm{V})
$$
by relative virtual localization, and derive the following identity:
\begin{proposition}
\begin{eqnarray*}
\textup{(\ref{eqn:KG})}&&
K^\bu_\mm(\lam)=\sum_{|\nu^\pm|=|\mu^\pm|}\Phi^\bu_{\mu^+,\nu^+}(-\sqrt{-1}\tau\lam)
z_{\nu^+} G^\bu_\nn(\lam;\tau)z_{\nu^-}\Phi^\bu_{\nu^-,\mu^-}(\frac{-\sqrt{-1}}{\tau}\lam).
\end{eqnarray*}
\end{proposition}

\subsection{Torus fixed points}

Given a morphism
$$
\Cxy\to X[m^+,m^-]
$$
which represents a point in $\MX^T$, let
$$
\tf=\pi[m^+,m^-]\circ f: C\to X.
$$
Then
\begin{eqnarray*}
&& \tf(C)\subset D_1\cup D_2\cup D_3\cup D_4\cup D_5,\\
&&\tf(x_i)\in\{ z_+,\tz_+\},\ \ \ \tf(y_j)\in\{ z_-,\tz_-\}.
\end{eqnarray*}
See Figure 3 in Section \ref{action} for the configuration of the
$T$-invariant divisors $D_1,\ldots,D_5$ and the $T$ fixed points
$z_0,z_+,z_-,\tz_+,\tz_-$.

If
$$
\tf_*(C)=n_1D_1+ n_2 D_2 + n_3 D_3 + n_4 D_4 + n_5 D_5
$$
as divisors, then
$$
\tf_*[C]=(n_1+n_2)\beta_1 + (n_3+n_4)\beta_3 + (n_5-n_2-n_4)\beta_5
$$
as homology classes.

Let
$$
J=\{\five{n}\in\bZ^5\mid n_i\geq 0, n_1+n_2=|\mu^-|, n_3+n_4=|\mu^+|,
n_5=n_2+n_4\}.
$$
Given $\hat{n}=\five{n}\in J$, let
$$
\cM_{\hat{n}}\subset \MX^T
$$
be the subset which corresponds to
$$
\tf_*(C)=n_1D_1+ n_2 D_2 + n_3 D_3 + n_4 D_4 + n_5 D_5.
$$

Then $\MX^T$ is a disjoint union of
$$
\{\cM_{\hat{n}}:\hat{n}\in J\}.
$$
We have the following vanishing lemma:
\begin{lemma}\label{vanish}
Let $\hat{n}\in J$, and let
$i_{\hat{n}}:\cM_{\hat{n}}\to \MX^T$ be the inclusion. Then
$$
i^*_{\hat{n}}e_T(\xm{V})=0
$$
unless $\hat{n}=(|\mu^-|,0,|\mu^+|,0,0)$.
\end{lemma}
\begin{proof}

We use the notation in Section \ref{obstruction}. Let
$L=\cO_X(-D_1-D_3)$. We have the following short
exact sequence of sheaves on $\MX$:
\begin{equation}\label{eqn:split}
0\to \tF^*L(-\xm{\cD})\to \tF^* L \to (\tF^*L)_{\xm{\cD}}\to 0.
\end{equation}

Let
$$
s_i:\MX\to\xm{\cU}
$$
be the section corresponds to the $i$-th marked point,
$$
\mathrm{ev}_i=\tF\circ s_i:\MX\to X
$$
be evaluation at the $i$-th marked point.
Then (\ref{eqn:split}) gives the following long exact sequence:
\begin{eqnarray*}
0&\to& R^0\pi_*\tF^*L(-\xm{D})\to R^0\tF^* L \to
\bigoplus_{i=1}^{l(\mu^+)+l(\mu^-)}\mathrm{ev}_i^* L\\
&\to& R^1\pi_*\tF^*L(-\xm{D})\to R^1\tF^* L \to 0.
\end{eqnarray*}
We have $R^0\tF^*L=0$, so
$$
\xm{\tV}=R^1\tF^*L
$$
is a vector bundle over $\MX$. We have the following short exact sequence of
vector bundles over $\MX$:
\begin{equation}\label{eqn:vtv}
0\to \bigoplus_{i=1}^{l(\mu^+)+l(\mu^-)}\mathrm{ev}_i^* L \to \xm{V}\to
\xm{\tV}\to 0.
\end{equation}

The restriction of the above exact sequence to $\cM_{\hat{n}}$ is
$$
0\to L_{z_-}^{\oplus l(\si^1)}\oplus L_{\tz_-}^{\oplus l(\si^2)}
    L_{z_+}^{\oplus l(\si^3)}\oplus L_{\tz_-}^{\oplus l(\si^4)}
   \to i^*_{\hat{n}}\xm{V}\to i^*_{\hat{n}}\xm{\tV}\to 0,
$$
where $\si^1,\ldots,\si^4$ are partitions  determined by
\begin{equation}\label{sigma}
\{\mu_j^-:\tf(y_j)\in z_-\},\ \ \
\{\mu_j^-:\tf(y_j)\in \tz_-\},\ \ \
\{\mu_i^+:\tf(x_i)\in z_+\},\ \ \
\{\mu_i^+:\tf(x_i)\in \tz_+\},
\end{equation}
respectively. Note that $\si^1,\si^2,\si^3,\si^4$ are constant on each connected components
of $\cM_{\hat{n}}$, and
$$
\si^1\cup\si^2=\mu^-,\ \ \, \si^3\cup\si^4=\mu^+.
$$

We have seen in Section \ref{action} that
$$
e_T(L_{z_+})=-\beta,\ \ \ e_T(L_{z_-})=-\alpha, \ \ \
e_T(L_{\tz^+})=e_T(L_{\tz^-})=0,
$$
so
$$
i^* e_T(\xm{V})=0
$$
unless 
\begin{equation} \label{intersect}
(\si^1,\si^2,\si^3,\si^4)=(\mu^-,\emptyset,\mu^+,\emptyset).
\end{equation}
 
Let $\hat{n}=\five{n}\in J$, $n_5\neq 0$.
Let $\cM_{\hat{n}}(k)$ be the subset of $\cM_{\hat{n}}$ which consists of
points
$$
\left[ f:C\to X[m^+,m^-]\right]\in \cM_{\hat{n}}
$$
such that (\ref{intersect}) is true, and
$$
f^{-1}(D_5-\{\tz^+,\tz^-\})
$$
has $k$ connected components, where $1\leq k\leq n_5$.
Each $\cM_{\hat{n}}(k)$ is a union of connected components
of $\cM_{\hat{n}}$.

We claim that 
$$
e_T(\xm{V})|_{\cM_{\hat{n}}(k)}=0
$$
for all $\hat{n}=\five{n}\in J$, $n_5\neq 0$, $k=1,\ldots,n_5$.
This will complete the proof.

Let
$$
\left[f:C\to X[m^+,m^-]\right]\in\cM_{\hat{n}}(k).
$$
Then $C=C_1\cup C_2$, where $C_1$ is the closure of
$f^{-1}(D_5-\{\tz^+,\tz^-\})$, which is a disjoint
union of $k$ projective lines, and $C_2$ is the union of 
other irreducible components of $C$. By (\ref{intersect}), 
the ramification divisor $R\subset C_2$, and
$C_1$ and $C_2$ intersect at $2k$ nodes. We have
\begin{eqnarray*}
0&\to& H^0(C,\tf^*L(-R))\to H^0(C_1,\tf^*L|_{C_1})\oplus H^0(C_2,\tf^* L(-R)|_{C_2}) \to
L_{\tz^+}^{\oplus k}\oplus L_{\tz^-}^{\oplus k}\\
&\to& H^1(C,\tf^*L(-R)) \to H^1(C_1,\tf^*L|_{C_1}) \oplus H^1(C_2,\tf^*L(-R)|_{C_2})\to 0
\end{eqnarray*}
where
$$
H^0(C,\tf^* L(-R))=0=H^0(C_2,\tf^*L(-R)|_{C_2}).
$$
The restriction of $L$ to $D_5$ is (equivariantly) trivial, so 
$$
H^0(C_1,\tf^*L|C_1)\cong L_{\tz^-}^{\oplus k},\ \
H^1(C_1,\tf^*L|C_1)=0.
$$
We have
$$
0\to L_{\tz^+}^{\oplus k}\to H^1(C,\tf^*L)\to H^1(C_2,\tf^*L)\to 0,
$$
so
$$
\xm{V}|_{\cM_{\hat{n}}(k)}= L_{\tz^+}^k\oplus V',
$$
and
$$
e_T(\xm{V})|_{\cM_{\hat{n}}(k)}=0.
$$
\end{proof}

Lemma \ref{vanish} tells us that $\cM_{\hat{n}}$ does not contribute to
the localization calculation of
$$
\xm{K}=\frac{1}{\amm}\int_{[\MX]^\vir}e_T(\xm{V})
$$
if $\hat{n}\neq (|\mu^-|,0,|\mu^+|,0,0)$.

\subsection{Admissible labels}

From now on, we only consider
$$
\hat{\cM}=\cM_{(|\mu^-|,0,|\mu^+|,0,0)}\subset \MX^T.
$$
Given a morphism
$$
\domain\to X[m^+,m^-]
$$
which represents a point in $\hat{\cM}$,
let
$$
C^0=\tf^{-1}(z_0), \ \ \
C^\pm=\tf^{-1}(z_\pm),
$$
where $z_0, z_+, z_-$ are defined as in Section \ref{action}.
Then
$$
C=C^+\cup L^+ \cup C^0\cup L^- \cup C^-,
$$
where $L^+, L^-$ are unions of projective lines,
$f|_{L^+}:L^+\to D_3$ is a degree $d^+=|\mu^+|$ cover
fully ramified over $z_0$ and $z_+$, and
$f|_{L^-}:L^-\to D_1$ is a degree $d^-=|\mu^-|$ cover
fully ramified over $z_0$ and $z_-$.

Define
$$
\bP^\pm(m^\pm)=\pi[m^+,m^-]^{-1}(z_\pm).
$$
Let
\begin{eqnarray*}
f^\pm = f|_{C^\pm}:C^\pm &\to& \bP^\pm(m^\pm),\\
\tf^+ = f|_{L^+}:L^+ &\to& D_3,\\
\tf^- = f|_{L^-}:L^- &\to& D_1.
\end{eqnarray*}
The degrees of $\tf^\pm$ restricted to irreducible components
of $L^\pm$ determine a partition $\nu^\pm$ of $d^\pm$.

Let $C_1^0,\ldots, C_k^0$ be the connected components of $C^0$, and
let $g_i$ be the arithmetic genus of $C_i^0$. (We define
$g_i=0$ if $C_i^0$ is a point.) Define
$$
\chi^0=\sum_{i=1}^k(2-2g_i).
$$

We define $\chi^+, \chi^-$ similarly. Then
$$
-\chi^+ + 2l(\nu^+) - \chi^0 + 2l(\nu^-) -\chi^- = -\chi.
$$

Note that $\chi^\pm \leq 2\min \{l(\mu^\pm),l(\nu^\pm)\}$.
So
$$
\nmp\geq 0,
$$
and the equality holds if and only if $m^+=0$. In this case, we have
$\nu^+=\mu^+$, $\chi^+= 2l(\mu^+)$. Similarly,
$$
\nmm\geq 0,
$$
and the equality holds if and only if $m^-=0$. In this case, we have
$\nu^- = \mu^-$, $\chi^-= 2l(\mu^-)$.
There are four cases:

\paragraph{Case 1: $m^+=m^-=0$.} Then $f^+, f^-$ are constant maps, and
$\nu^\pm=\mu^\pm$.

\paragraph{Case 2: $m^+>0, m^-=0$.} Then $f^-$ is a constant map,
$\nu^-=\mu^-$, and $f^+$ represents a point in
$$
\Mp//\bC^*
$$
up to an element in $\Aut(\nu^+)$.

\paragraph{Case 3: $m^+=0, m^->0$.} Then $f^+$ is a constant map,
$\nu^+ =\mu^+$, and $f^-$ represents a point in
$$
\Mm//\bC^*
$$
up to an element in $\Aut(\nu^-)$.

\paragraph{Case 4: $m^+, m^- >0$.} Then $f^+$ represents
a point in
$$
\Mp//\bC^*
$$
up to an element of $\Aut(\nu^+)$,
and $f^-$ represents an point in
$$
\Mm//\bC^*
$$
up to an element in $\Aut(\nu^-)$.

\begin{df}
An \emph{admissible label}  is a 5-uple $\lab$
such that
\begin{itemize}
\item $\chi^+,\chi^0,\chi^-\in 2\bZ$.
\item $\nu^\pm$ is a partition of $d^\pm$.
\item $\chi^0\leq 2\min\{l(\nu^+),l(\nu^-)\}$,
      $\chi^\pm\leq 2\min\{l(\mu^\pm), l(\nu^\pm)\}$.
\item $-\chi^+ +2 l(\nu^+) -\chi^0 +2 l(\nu^-) -\chi^- =  -\chi$.
\end{itemize}
Let $\GX$ denote the set of all admissible labels.
\end{df}

For a nonnegative integer $g$ and a positive integer $h$,
let $\Mbar_{g,h}$ be the moduli space of
stable curves of genus $g$ with $h$ marked points.
$\Mbar_{g,h}$ is empty for $(g,h)=(0,1), (0,2)$,
but we will assume that $\Mbar_{0,1}$ and $\Mbar_{0,2}$
exist and satisfy
\begin{eqnarray*}
\int_{\Mbar_{0,1}}\frac{1}{1-d\psi}&=&\frac{1}{d^2}\\
\int_{\Mbar_{0,2}}\frac{1}{(1-\mu_1\psi_1)(1
-\mu_2\psi_2)}
&=& \frac{1}{\mu_1+\mu_2}
\end{eqnarray*}
for simplicity of notation. Such an assumption will give
the correct final results.

For a nonnegative integer $g$ and a positive integer $h$,
let $\Mbar^\bullet_{\chi,h}$ be the moduli of possibly
disconnected  stable curves $C$  with $h$ marked points such that
\begin{itemize}
\item If $C_1,\ldots,C_k$ are connected components of $C$, and $g_i$ is the
  arithmetic genus of $C_i$, then
$$
\sum_{i=1}^k (2-2g_i)=\chi.
$$
\item Each connected component contains at least one marked point.
\end{itemize}

The connected components of $\Mbar^\bu_{\chi,h}$ are of the form
$$
\Mbar_{g_1,h_1}\times\cdots\times\Mbar_{g_k,h_k}.
$$
where
$$
\sum_{i=1}^k(2-2g_i)=\chi,\ \ \ \sum_{i=1}^k h_i=h.
$$
The restriction of the Hodge bundle $\bE\to \Mbar^\bu_{\chi,h}$ to
the above connected component is the direct sum of the Hodge bundles
on each factor, and
$$
\Lambda^\vee(u)=\prod_{i=1}^k\Lambda_{g_i}^\vee(u).
$$

We define
$$
\Mbar_{2l(\mu^+),\mu^+,\chi,\mu^-,2l(\mu^-)}=
\Mbar^\bu_{\chi,l(\mu^+)+l(\mu^-)}.
$$
For $-\chi^\pm+l(\mu^\pm)+l(\nu^\pm)>0$, we define
\begin{eqnarray*}
\Mbar_{\chi^+,\nu^+,\chi^0,\mu^-,2l(\mu^-)} &=&
\left(\Mp//\bC^*\right) \times \Mbar_{\chi^0,l(\nu^+)+l(\mu^-)},\\
\Mbar_{2l(\mu^+),\mu^+,\chi^0,\nu^-,\chi^-} &=&
\Mbar_{\chi^0,l(\mu^+)+l(\nu^-)}\times \left(\Mm//\bC^*\right),
\end{eqnarray*}
\begin{eqnarray*}
&&\xnxnx{\Mbar}\\
&=&\left(\Mp//\bC^*\right)
   \times\Mbar^\bullet_{\chi^0,l(\nu^+)+l(\nu^-)}
   \times \left(\Mm//\bC^*\right).
\end{eqnarray*}

For every
$\lab \in \GX$, there is a morphism
$$
\xnxnx{i}:\xnxnx{\Mbar}\to \MX,
$$
whose image $\xnxnx{\cF}$ is a union of connected components of $\hat{\cM}$.
The morphism $\xnxnx{i}$ induces an isomorphism
$$
\xnxnx{\Mbar}/\xnxnx{A} \cong \xnxnx{\cF},
$$
where $A_{2l(\mu^+),\mu^+,\chi,\mu^-,l(\mu^-)}$ is trivial, and
for $\chi^\pm+l(\mu^\pm)+l(\nu^\pm)>0$, we have
\begin{eqnarray*}
&& 1\to \prod_{i=1}^{l(\nu^+)}\bZ_{\nu_i^+}
\to A_{\chi^+,\nu^+,\chi^0,\mu^-,2l(\mu^-)}\to \Aut(\nu^+) \to 1,\\
&& 1\to  \prod_{j=1}^{l(\nu^-)}\bZ_{\nu_j^-}
\to A_{2l(\mu^+),\mu^+,\chi^0,\nu^-,\chi^-} \to  \Aut(\nu^-) \to 1,
\end{eqnarray*}
$$
1\to \prod_{i=1}^{l(\nu^+)}\bZ_{\nu_i^+}\times \prod_{j=1}^{l(\nu^-)}\bZ_{\nu_j^-}
\to \xnxnx{A}\to \Aut(\nu^+)\times \Aut(\nu^-) \to 1.
$$
So for $-\chi^\pm+l(\mu^\pm)+l(\nu^\pm)>0$, we have
$$
|A_{\chi^+,\nu^+,\chi^0,\mu^-,2l(\mu^-)}| = z_{\nu^+},\ \ \
|A_{2l(\mu^+),\mu^+,\chi^0,\nu^-,\chi^-}| = z_{\nu^-},\ \ \
$$
$$
|\xnxnx{A}| = z_{\nu^+}z_{\nu^-}.
$$

The stack $\hat{\cM}$ is a disjoint union of
$$
\{\xnxnx{\cF}:\lab\in \GX\}.
$$

\subsection{Contribution from each admissible label}
Let $\xnxnx{N}^\vir$ on $\xnxnx{\Mbar}$ be the pull back of the virtual normal
bundle of $\xnxnx{\cF}$ in $\MX$.
Calculations similar to those in \cite[Appendix A]{LLZ} show
that
$$
\frac{\xnxnx{i}^* e_T(V^\bu_{\chi,\mu^+,\mu^-}) }
{e_T(\xnxnx{N}^\vir)}
= A^+ A^0  A^-,
$$
where
\begin{eqnarray*}
A^0&=&(-1)^{|\nu^+|+|\nu^-|+1} a_{\nu^+} a_{\nu^-}\prod_{i=1}^{l(\nu^+)}
\frac{\prod_{a=1}^{\nu^+_i-1}(\nu^+_i\beta + a\alpha)}
     {(\nu^+_i-1)!\alpha^{\nu^+_i-1} }
\prod_{j=1}^{l(\nu^-)}
\frac{\prod_{a=1}^{\nu^-_j-1}(\nu^-_j\alpha + a\beta)}
     {(\nu^-_j-1)!\beta^{\nu^-_j-1} }\\
&&\cdot \frac{\Lambda^\vee(\alpha)\Lambda^\vee(\beta)
\Lambda^\vee(-\alpha-\beta)(\alpha\beta(\alpha+\beta))^{l(\nu^+)+l(\nu^-)-1} }
{\prod_{i=1}^{l(\nu^+)}\alpha(\alpha-\nu^+_i\psi_i)
 \prod_{j=1}^{l(\nu^-)}\beta(\beta-\nu^-_j\psi_{l(\nu^+)+j})}\\
A^+&=& \left\{\begin{array}{ll}(-1)^{l(\mu^+)}, & \chi^-=2l(\mu^-)\\
(-1)^{-\frac{\chi^+}{2}+l(\nu^+)+l(\mu^+)}a_{\nu^+}
\frac{\beta^{-\chi^++l(\mu^+)+l(\nu^+)}}{-\alpha-\psi^+}, &\textup{otherwise}
\end{array}\right.\\
A^-&=& \left\{\begin{array}{ll}(-1)^{l(\mu^-)},& \chi^+=2l(\mu^+)\\
(-1)^{-\frac{\chi^-}{2}+l(\nu^-)+l(\mu^-)}a_{\nu^-}
\frac{\alpha^{-\chi^-+l(\mu^-)+l(\nu^-)}}{-\beta-\psi^-}, &\textup{otherwise}
\end{array}\right.
\end{eqnarray*}

From the definitions in Section \ref{sec:hodge} and
Proposition \ref{pro:xELSV}, we have
\begin{eqnarray*}
& &\xnxnx{I}(\alpha,\beta)\\
&=&\frac{1}{|\xnxnx{A}|}
   \int_{[\xnxnx{\Mbar}]^\vir}
\frac{\xnxnx{i}^* e_T(V^\bu_{\chi,\mu^+,\mu^-})}{e_T(\xnxnx{N}^\vir)}\\
&=&\amm\frac{ \sqrt{-1}^{l(\mu^+)+l(\mu^-)} }{ (-1)^{|\mu^+|+|\mu^-|} }
G^\bu_{\chi^0,\nu^+,\nu^-}(\alpha,\beta)\\
&& \cdot z_{\nu^+} \frac{(-\sqrt{-1}\beta/\alpha)^{-\chi^+ +l(\nu^+) +l(\mu^+)} }
{(-\chi^+ +l(\nu^+) +l(\mu^+))!}H^\bu_{\chi^+,\nu^+,\mu^+}
\cdot  z_{\nu^-}
\frac{(-\sqrt{-1}\alpha/\beta)^{-\chi^- +l(\nu^-) +l(\mu^-)}}
{(-\chi^- +l(\nu^-) +l(\mu^-))!}H^\bu_{\chi^-,\nu^-,\mu^-}
\end{eqnarray*}
Let $\tau=\beta/\alpha$. Then
\begin{eqnarray*}
& &\xnxnx{I}(\alpha,\beta)\\
&=&\amm\frac{\sqrt{-1}^{l(\mu^+)+l(\mu^-)} }{(-1)^{|\mu^+|+|\mu^-|}}
G^\bu_{\chi^0,\nu^+,\nu^-}(\tau)\\
&&\cdot z_{\nu^+} \frac{(-\sqrt{-1}\tau)^{-\chi^+ +l(\nu^+) +l(\mu^+)} }
{(-\chi^+ +l(\nu^+) +l(\mu^+))!}H^\bu_{\chi^+,\nu^+,\mu^+}
 \cdot z_{\nu^-} \frac{(-\sqrt{-1}/\tau)^{-\chi^- +l(\nu^-) +l(\mu^-)}}
{(-\chi^- +l(\nu^-) +l(\mu^-))!}H^\bu_{\chi^-,\nu^-,\mu^-}
\end{eqnarray*}

\subsection{Sum over admissible labels}

\begin{eqnarray*}
& & \xm{K}\\
&=&\frac{1}{\amm}\int_{[\MX]^\vir}e_T(\xm{V})\\
&=&\frac{1}{\amm}\sum_{\lab\in \GX}\frac{1}{|\xnxnx{A}|}\\
&&\cdot \int_{[\xnxnx{\Mbar]}^\vir}
\frac{\xnxnx{i}^* e_T(V^\bu_{\chi,\mu^+,\mu^-}) }
{e_T(\xnxnx{N}^\vir)}\\
&=&\frac{1}{\amm}\sum_{\lab\in \GX}\xnxnx{I}(\alpha,\beta)\\
&=&\frac{\sqrt{-1}^{l(\mu^+)+l(\mu^-)} }{(-1)^{|\mu^+|+|\mu^-|} }
\left(\sum_{\lab\in\GX}G^\bu_{\chi^0,\nu^+,\nu^-}(\tau)\right.\\
&&\left.\cdot z_{\nu^+}
\frac{(-\sqrt{-1}\tau)^{-\chi^+ +l(\nu^+) +l(\mu^+)}}
{(-\chi^+ +l(\nu^+) +l(\mu^+))!}H^\bu_{\chi^+,\nu^+,\mu^+}
\cdot z_{\nu^-}
\frac{(-\sqrt{-1}/\tau)^{-\chi^- +l(\nu^-) +l(\mu^-)}}
{(-\chi^- +l(\nu^-) +l(\mu^-))!}H^\bu_{\chi^-,\nu^-,\mu^-}\right)
\end{eqnarray*}

Recall that
$$
K^\bu_\mm(\lambda)=
\sum_{\chi\in 2\bZ,\chi\leq 2(l(\mu^+)+l(\mu^-))}\lam^{\dis}
\frac{(-1)^{|\mu^+|+|\mu^-|}}{\sqrt{-1}^{l(\mu^+)+l(\mu^-)}}\xm{K}.
$$
We have
$$
K^\bu_\mm(\lambda)=\sum_{|\nu^\pm|=|\mu^\pm|}
\Phi^\bu_{\mu^+,\nu^+}(-\sqrt{-1}\tau\lam)
z_{\nu^+}G^\bu_{\nn}(\lambda;\tau)z_{\nu^-}
\Phi^\bu_{\nu^-,\mu^-}\left(\frac{-\sqrt{-1}}{\tau}\lam\right).
$$
This finishes the proof of (\ref{eqn:KG}).


\begin{thebibliography}{99}

\bibitem{AKMV}
M.~Aganagic, A.~Klemm, M.~Marino, C.~Vafa,
{\em The topological vertex},
preprint, hep-th/0305132.
 
\bibitem{Aga-Mar-Vaf}
M.~Aganagic, M.~Marino, C.~Vafa,
{\em All loop topological string amplitudes from Chern-Simons theory},
preprint, hep-th/0206164.

\bibitem{Bry-Pan} J.~Bryan, R.~Pandharipande,
{\em Curves in Calabi-Yau 3-folds and topological quantum field theory},
preprint, math.AG/0306316.

\bibitem{Dia-Flo}
D.-E.~Diaconescu, B.~Florea,
{\em Localization and gluing of topological amplitudes},
preprint, hep-th/0309143.

 
\bibitem{Dij}
R.~Dijkgraaf,
Mirror symmetry and elliptic curves, {\em The moduli space of curves},
R.~Dijkgraaf, G.~van der Geer (editors), Prog. in Math., {\bf 129}, Birkha\"user, 1995.



\bibitem{ELSV}
T.~Ekedahl, S.~Lando, M.~Shapiro, A.~Vainshtein,
{\em Hurwitz numbers and intersections on moduli spaces of curves}.
{Invent. Math}. {\bf 146} (2001), no. 2, 297--327.

\bibitem{Fan-Pan}
B.~Fantecchi, R.~Pandharipande,
{\em Stable maps and branch divisors},
Compositio Math. {\bf 130} (2002), no. 3, 345--364.


\bibitem{Gou-Jac-Vai}
I.P.~Goulden, D.M.~Jackson, A.~Vainshtein,
{\em The number of ramified coverings of the sphere by the torus and surfaces of higher genera},
Ann. of Comb. {\bf 4} (2000), 27-46.

\bibitem{Gra-Pan}
T.~Graber, R.~Pandharipande,
{\em Localization of virtual classes}.
Invent. Math. {\bf 135} (1999), no. 2, 487--518.

\bibitem{Gra-Vak1}
T.~Graber, R.~Vakil,
{\em Hodge integrals and Hurwitz numbers via virtual localization}.
Compositio Math. {\bf 135} (2003), no. 1, 25--36.

\bibitem{Gra-Vak2}
T.~Graber, R.~Vakil,
{\em Relative virtual localization and vanishing of tautological
classes on moduli spaces of curves}, preprint, math.AG/0309227.


\bibitem{Iqb}
A.~Iqbal,
{\em All genus topological amplitudes and $5$-brane webs as Feynman diagrams},
preprint, hep-th/0207114.

\bibitem{Ion-Par1}E.-N. Ionel, T. Parker,
{\em Relative Gromov-Witten invariants},
Ann. of Math. (2) {\bf 157} (2003), no. 1, 45--96. 
 
\bibitem{Ion-Par2}
E.-N. Ionel, T. Parker,
{\em The symplectic sum formula for Gromov-Witten invariants},
preprint, math.SG/0010207.

\bibitem{Kat-Liu}
S.~Katz, C.-C.~Liu,
{\em Enumerative geometry of stable maps with Lagrangian boundary conditions and multiple covers
of the disc},
Adv. Theor. Math. Phys. {\bf 5} (2001), 1-49.

\bibitem{Kon}
M.~Kontsevich,
{\em Intersection theory on the moduli space of curves and the matrix Airy function}.
Comm. Math. Phys. {\bf 147} (1992), no. 1, 1--23.

\bibitem{Li-Rua} A.M.~Li, Y.~Ruan,
{\em Symplectic surgery and Gromov-Witten invariants of Calabi-Yau
  3-folds},
Invent. Math.{\bf 145} (2001), no. 1, 151--218. 

\bibitem{Li-Zha-Zhe}
A.M.~Li, G.~Zhao, Q.~Zheng,
{\em The number of ramified coverings of a Riemann surface by Riemann surface},
Comm. Math. Phys. {\bf 213} (2000), no. 3, 685--696.

\bibitem{Li1}
J.~Li,
{\em Stable Morphisms to singular schemes and relative stable morphisms},
J. Diff. Geom. {\bf 57} (2001), 509-578.

\bibitem{Li2}
J.~Li,
{\em Relative Gromov-Witten invariants and a degeneration formula of Gromov-Witten invariants},
J. Diff. Geom. {\bf 60} (2002), 199-293.

\bibitem{Li3}
J.~Li,
{\em Lecture notes on relative GW-invariants},
preprint.


\bibitem{LLLZ}
J.~Li, C.-C.~Liu, K.~Liu, J.~Zhou,
{\em A mathematical theory of the topological vertex}, in preparation.

\bibitem{LLZ}
C.-C.~Liu, K.~Liu, J.~Zhou,
{\em A proof of a conjecture of Mari\~no-Vafa on Hodge Integrals},
preprint, math.AG/0306434.


\bibitem{Mac}
I.G.~MacDonald,
{\em Symmetric functions and Hall polynomials}, 2nd edition.Claredon Press, 1995.

\bibitem{Mar-Vaf}
M.~Mari\~{n}o, C.~Vafa,
{\em Framed knots at large $N$},
Orbifolds in mathematics and physics (Madison, WI, 2001),
185--204, Contemp. Math., 310, Amer. Math. Soc., Providence, RI, 2002.

\bibitem{Mor-Luk}
H.R.~Morton, S.G.~Lukac,
{\em The HOMFLY polynomial of the decorated Hopf link},
preprint, math.GT/0108011.


\bibitem{Oko-Pan}
A.~Okounkov, R.~Pandharipande,
{\em Hodge integrals and invariants of the unknots},
preprint, math.AG/0307209.

\bibitem{Wit1}
E.~Witten, {\em Quantum field theory and the Jones polynomial},
Commun. Math. Phys. {\bf 121} (1989) 351-399.


\bibitem{Wit}
E.~Witten,
{\em Two-dimensional gravity and intersection theory on moduli space}.
Surveys in differential geometry (Cambridge, MA,
1990), 243--310, Lehigh Univ., Bethlehem, PA, 1991.


\bibitem{Zho1}
J.~Zhou,
{\em Hodge integrals, Hurwitz numbers, and symmetric groups},
preprint, math.AG/0308024.

\bibitem{Zho2}
J.~Zhou
{\em Marino-Vafa formula and BPS numbers in local $\bP^2$ and $\bP^1 \times \bP^1$ geometry},
preprint, submitted.

\bibitem{Zho3}
J.~Zhou,
{\em A conjecture on Hodge integrals},
preprint, math.AG/0310282.

\bibitem{Zho4}
J.~Zhou,
{\em Localizations on moduli spaces and free field realizations of Feynman rules},
preprint, math.AG/0310283.

\end{thebibliography}
\end{document}